\documentclass{article}
\usepackage[utf8]{inputenc}

\usepackage{amssymb}
\usepackage{amsmath}
\usepackage{amsthm}
\usepackage{amsfonts}
\usepackage{tikz}
\usepackage{xcolor}
\usepackage[hidelinks]{hyperref}
\usepackage{cleveref}
\usepackage{titling}
\thanksmarkseries{arabic}

\usepackage[top=1.0in, bottom=1.0in, left=1.0in, right=1.0in, headheight=1.0in]{geometry}
\usepackage{mathtools}
\usepackage{youngtab}
\usepackage{graphicx}
\graphicspath{ {./Graphs/} }

\usepackage[vlined]{algorithm2e}
\usepackage{enumerate}
\usepackage[normalem]{ulem}
\usepackage{multirow}

\usepackage{soul}

\usepackage{tikz}
\usetikzlibrary{arrows.meta, calc}

\usepackage{subcaption}

\theoremstyle{plain}
\newtheorem{theorem}{Theorem}
\newtheorem{lemma}[theorem]{Lemma}
\newtheorem{proposition}[theorem]{Proposition}
\newtheorem{conjecture}[theorem]{Conjecture}

\newtheorem{corollary}[theorem]{Corollary}

\newtheorem{question}[theorem]{Question}

\theoremstyle{definition}

\DeclareMathOperator{\diam}{diam}
\DeclareMathOperator{\ecc}{ecc}

\newcommand{\x}{\times}

\DeclareMathOperator{\rad}{rad}
\newcommand{\MaxTrail}{\ell}

\usepackage{xcolor}
\definecolor{nicepurple}{RGB}{200, 90, 160}

\def\pone{{\rm{P}\ensuremath{1}}}
\def\ptwo{{\rm{P}\ensuremath{2}}}
\def\poneposs{{\rm{P}\ensuremath{1}'s}}
\def\ptwoposs{{\rm{P}\ensuremath{2}'s}}

\title{Trail Trap: a variant of Partizan Edge Geography}
\author{Calum Buchanan\thanks{Department of Mathematics \& Statistics, University of Vermont, Burlington, VT, USA {\tt calum.buchanan@uvm.edu}}\and MacKenzie Carr\thanks{Department of Mathematics, Simon Fraser University, Burnaby, Canada {\tt mackenzie\_carr@sfu.ca}}\and Alexander Clifton\thanks{Department of Theoretical Computer Science, Czech Technical University in Prague, Prague, Czechia {\tt alexander.clifton@fit.cvut.cz}}\and Stephen G. Hartke\thanks{Department of Mathematical and Statistical Sciences, University of Colorado Denver, USA \texttt{stephen.hartke@ucdenver.edu}}\and Vesna Ir\v{s}i\v{c} Chenoweth\thanks{University of Ljubljana, Slovenia \tt vesna.irsic@fmf.uni-lj.si}\and Nicholas Sieger\thanks{University of California, San Diego, CA, USA \tt nsieger@ucsd.edu}\and Rebecca Whitman\thanks{University of California, Berkeley, CA, USA \tt rebecca\_whitman@berkeley.edu}}
\date{\today}

\begin{document}
\maketitle

\begin{abstract}
    We study a two-player game played on undirected graphs called {\sc Trail Trap}, which is a variant of a game known as {\sc Partizan Edge Geography}. One player starts by choosing any edge and moving a token from one endpoint to the other; the other player then chooses a different edge and does the same. Alternating turns, each player moves their token along an unused edge from its current vertex to an adjacent vertex, until one player cannot move and loses. We present an algorithm to determine which player has a winning strategy when the graph is a tree and partially characterize the trees on which a given player wins. Additionally, we show that it is NP-hard to determine if Player~2 has a winning strategy on {\sc Trail Trap} from the starting position, even for connected bipartite planar graphs with maximum degree $4$. We determine which player has a winning strategy for certain subclasses of complete bipartite graphs and grid graphs, and we propose several open problems for further study. 
\end{abstract}

\section{Introduction}

Let $G = (V(G),E(G))$ be a graph (all graphs we consider are undirected, finite, and simple). A \emph{walk} on $G$ is a sequence $v_1 v_2 \ldots v_{m}$ of vertices in $V(G)$ such that $v_i v_{i+1} \in E(G)$ for all $1 \le i \le m - 1$. A \emph{trail} is a walk with no repeated edges and a \emph{path} is a trail with no repeated vertices. In this work, we study a perfect information game, {\sc Trail Trap}, first introduced in \cite{ArCaDa21}, in which two players compete to build edge-disjoint trails on $G$. 

To begin a game of {\sc Trail Trap} on a graph $G$, Player~1 (\pone{}) chooses a starting vertex on which to place their token, then moves the token along an incident edge $e$ (should such an edge exist) to an adjacent vertex, completing their turn. Player~2 (\ptwo{}) then places their token on any vertex, and moves it along an incident edge $f \not = e$ (should such an edge exist), to its second endpoint. Continuing with \pone{}, the players alternate turns, each moving their token from its current vertex along an unused edge to a neighboring vertex. A player cannot make a move along an edge previously used by either player, but there are no restrictions on vertices: a player may place their token on the same vertex multiple times throughout the game, and the two tokens may simultaneously occupy a vertex. We apply the Normal Play convention: the first player unable to make a move loses the game. Since $|E(G)|$ is finite, {\sc Trail Trap} cannot end in a draw. 

We denote a \emph{move} on a graph $G$, the act of moving a token from one vertex $u$ to an adjacent vertex $v$, by $u\to v$. Interchangeably, we say that a player \emph{moves} from $u$ to $v$. An edge $uv$ is \emph{used} if either player has moved from $u$ to $v$ or from $v$ to $u$ on any previous turn, and \emph{unused} otherwise. Each player's sequence of moves, at every turn, forms a trail through the graph. Thus, rather than moving tokens, we can think of the players as constructing two edge-disjoint trails on the graph. 

We provide an example of gameplay on the ``diamond" graph $G$, illustrated in Figure \ref{fig:example_game_diamond}. Suppose \pone{} first makes the move $b \to c$, and suppose \ptwo{} begins their trail with $c \to d$. From here, \poneposs{} second move must be $c \to a$, since this is the only unused edge incident to $c$. Similarly, \ptwoposs{} second move is $d \to b$, and \pone{} then makes the move $a \to b$. Since \ptwo{} has no available moves out of vertex $b$ (indeed, all the edges in the graph are used), the game is over and \pone{} wins. 

\begin{figure}[ht]
    \centering
    \begin{tikzpicture}
    [very thick,every node/.style={circle, draw=black, fill=black, inner sep=2}]

    \node[label=below:d] (1a) at (0,0){};
    \node[label=left:b] (2a) at (-1,1.5){};
    \node[label=right:c] (3a) at (1,1.5){};
    \node[label=a] (4a) at (0,3){};

    \draw[-{Latex[length=3mm]},color=red] (2a) -- node[draw=none,fill=none,above]{\pone{}}(3a);
    \draw (3a) -- (1a);
    \draw (3a) -- (4a);
    \draw (1a) -- (2a);
    \draw (4a) -- (2a);

    \node[label=below:d] (1b) at (3.25,0){};
    \node[label=left:b] (2b) at (2.25,1.5){};
    \node[label=right:c] (3b) at (4.25,1.5){};
    \node[label=a] (4b) at (3.25,3){};

    \draw[-{Latex[length=3mm]},color=red] (2b) -- node[draw=none,fill=none,above]{\pone{}}(3b);
    \draw[-{Latex[length=3mm]},color=cyan!70!black] (3b) -- node[draw=none,fill=none,right]{\ptwo{}} (1b);
    \draw (3b) -- (4b);
    \draw (1b) -- (2b);
    \draw (4b) -- (2b);

    \node[label=below:d] (1c) at (6.5,0){};
    \node[label=left:b] (2c) at (5.5,1.5){};
    \node[label=right:c] (3c) at (7.5,1.5){};
    \node[label=a] (4c) at (6.5,3){};

    \draw[-{Latex[length=3mm]},color=red] (2c) -- node[draw=none,fill=none,above]{\pone{}}(3c);
    \draw[-{Latex[length=3mm]},color=cyan!70!black] (3c) -- node[draw=none,fill=none,right]{\ptwo{}} (1c);
    \draw[-{Latex[length=3mm]},color=red] (3c) -- node[draw=none,fill=none,right]{\pone{}} (4c);
    \draw (1c) -- (2c);
    \draw (4c) -- (2c);

    \node[label=below:d] (1d) at (9.75,0){};
    \node[label=left:b] (2d) at (8.75,1.5){};
    \node[label=right:c] (3d) at (10.75,1.5){};
    \node[label=a] (4d) at (9.75,3){};

    \draw[-{Latex[length=3mm]},color=red] (2d) -- node[draw=none,fill=none,above]{\pone{}}(3d);
    \draw[-{Latex[length=3mm]},color=cyan!70!black] (3d) -- node[draw=none,fill=none,right]{\ptwo{}} (1d);
    \draw[-{Latex[length=3mm]},color=red] (3d) -- node[draw=none,fill=none,right]{\pone{}} (4d);
    \draw[-{Latex[length=3mm]},color=cyan!70!black] (1d) -- node[draw=none,fill=none,left]{\ptwo{}} (2d);
    \draw (4d) -- (2d);

    \node[label=below:d] (1e) at (13,0){};
    \node[label=left:b] (2e) at (12,1.5){};
    \node[label=right:c] (3e) at (14,1.5){};
    \node[label=a] (4e) at (13,3){};

    \draw[-{Latex[length=3mm]},color=red] (2e) -- node[draw=none,fill=none,above]{\pone{}}(3e);
    \draw[-{Latex[length=3mm]},color=cyan!70!black] (3e) -- node[draw=none,fill=none,right]{\ptwo{}} (1e);
    \draw[-{Latex[length=3mm]},color=red] (3e) -- node[draw=none,fill=none,right]{\pone{}} (4e);
    \draw[-{Latex[length=3mm]},color=cyan!70!black] (1e) -- node[draw=none,fill=none,left]{\ptwo{}} (2e);
    \draw [-{Latex[length=3mm]},color=red] (4e) -- node[draw=none,fill=none,left]{\pone{}} (2e);

    \end{tikzpicture}
    \caption{Sample gameplay on the diamond graph, with five consecutive game states shown.}
    \label{fig:example_game_diamond}
\end{figure}

However, \ptwoposs{} first move $c \to d$ wasn't optimal, and with a different first move, \ptwo{} can ensure victory. In fact, with optimal play, \ptwo{} can win on $G$ no matter \poneposs{} first move. We leave this as an exercise to the reader. 


When a game of {\sc Trail Trap} ends, either \poneposs{} trail is of length one larger than \ptwoposs{} (and \pone{} is the winner), or the trails have the same length (and \ptwo{} is the winner).
In general, we say that a player has a \emph{winning strategy} on a graph $G$ if they can win {\sc Trail Trap} on $G$ no matter the other player's moves. If \pone{} has a winning strategy on $G$, then we say that $G$ is \emph{\pone{}-win} (\emph{\ptwo{}-win} is defined analogously). By the Fundamental Theorem of Combinatorial Games (see \cite{AlNoWo19}), exactly one of \pone{} and \ptwo{} has a winning strategy on any given graph. In this paper, we address the problem of characterizing \pone{}-win graphs for several families, including complete bipartite graphs and trees. We give a polynomial-time algorithm to determine the winner of {\sc Trail Trap} on an arbitrary tree. On the other hand, we prove that it is NP-hard to determine if \ptwo{} has a winning strategy even for connected planar bipartite graphs with maximum degree $4$. 

In the following subsection, we provide a brief survey of Geography games, of which {\sc Trail Trap} can be considered a variant. Then, in Section \ref{sec preliminaries}, we provide a necessary and sufficient condition for a disconnected graph to be \pone{}-win as well as a sufficient condition, using automorphisms, for a graph to be \ptwo{}-win.
In Section \ref{sec solved families}, we give exact win conditions for subsets of complete graphs, complete bipartite graphs, grid graphs, and prism graphs. In Section \ref{sec tree algorithm}, we characterize the \pone{}-win trees and give an $O(n^2)$ algorithm to determine the winner for trees of order $n$. In Section \ref{sec NP-hardness}, we prove that determining whether \ptwo{} wins is NP-hard  on a restricted class of graphs. Open questions arise in each section, and we summarize these in Section~\ref{section:conclusion}.

\subsection{Related Geography games}\label{sub:related games}

{\sc Trail Trap} can be seen as a variant of (Generalized) Geography, part of a larger family of walk-construction games (see~\cite{bodlaender1993complexity}). Geography is a two-player game in which players take turns moving a single, shared token from one vertex to an adjacent vertex, beginning at a fixed starting vertex (or root). In {\sc Vertex Geography}, no vertex can be repeated, and in {\sc Edge Geography}, no edge can be repeated. As in {\sc Trail Trap}, the winner is the last player able to make a legal move. Both {\sc Vertex Geography} and {\sc Edge Geography} can be played on directed or undirected graphs. 

{\sc Vertex Geography} is the simplest game in the Geography family; Fraenkel, Scheinerman, and Ullman~\cite{FSU93} provide an exact characterization of when each player wins on an undirected graph $G$ with a given starting vertex $v$. Specifically, they show that $(G,v)$ is \pone{}-win if and only if every maximum matching in $G$ includes an edge incident to $v$. In {\sc Unrooted Vertex Geography}, a generalization where \poneposs{} first move is to choose the starting vertex (and then \ptwo{} makes the first ``full" move to a new vertex), the graph $G$ is \ptwo{}-win if and only if it has a perfect matching. It follows that the winners of both of these variants can be computed in polynomial time \cite{edmonds65}. Monti and Sinaimeri~\cite{MS18} consider a modified version of {\sc Vertex Geography} where the goal is to determine whether \pone{} can force a win within $k$ moves. They show that this problem is PSPACE-complete, even for bipartite graphs, but can be solved in polynomial time for trees. 

{\sc Edge Geography} is, in general, more complicated than {\sc Vertex Geography.} The winner for {\sc Edge Geography} on undirected graphs can be computed in polynomial time for a few special families, including bipartite and complete graphs \cite{FSU93}. However, {\sc Undirected Edge Geography} in general is PSPACE-complete, even for subcubic planar graphs. The same holds for an unrooted generalization defined analogously to {\sc Unrooted Vertex Geography.} On directed graphs, both {\sc Vertex Geography} and {\sc Edge Geography} are PSPACE-complete (for these general results, see Lichtenstein and Sipser~\cite{LS80}, and Schaefer~\cite{Schaefer78}). 

In \emph{partizan} variants of Geography, players have distinct options for their moves, either because the players have different (and often complex) rulesets (for instance, as in \cite{Bosboom2020}) or because the players are using separate tokens. \emph{Partizan} games, where the two players' moves are no longer symmetric, are much more challenging to analyze than the impartial games mentioned above \cite{Conway00}. 
In~\cite{FS93}, Fraenkel and Simonson considered the (rooted) games {\sc Partizan Vertex Geography} and {\sc Partizan Edge Geography}, in which the two players are each assigned a starting vertex as well as their own token to move along unused edges in a graph without repeating a vertex (respectively, an edge).
They showed that the winning player of either game can be determined in $O(n^2)$ time for directed trees and $O(n)$ time for undirected trees, where $n$ is the order of the tree. However, no other classes are known to be solvable in polynomial time. Fox and Geissler~\cite{FG22} show that {\sc Partizan Vertex Geography} is PSPACE-complete for undirected bipartite graphs. {\sc Partizan Vertex Geography} and {\sc Partizan Edge Geography} are both NP-hard for directed acyclic graphs and directed planar graphs \cite{FS93}. Both games are known to be PSPACE-complete on bipartite digraphs with bounded degree: specifically, where the in-degree and out-degree of each vertex is at most two, and their sum is at most three \cite{FS93}.

{\sc Trail Trap} is formally a variant of {\sc Unrooted Partizan Edge Geography}. In this variation, each player places a token \emph{and} moves along an edge on their first turn. This contrasts with the unrooted (impartial) games mentioned in \cite{FSU93}, where \pone{} places the token on their first turn, and \ptwo{} is the first to move the token. The additional degree of freedom for the players in {\sc Trail Trap} can change the outcome of the game. For example, the path on three vertices is \pone{}-win in the game where both players place tokens before either moves, but \ptwo{}-win in {\sc Trail Trap}. {\sc Trail Trap} principally differs from these other games in its wide range of starting positions, and indeed much of our analysis (for instance, in \Cref{thm tree characterization}) focuses entirely on the players' first moves. Most other geography variants are PSPACE-complete, and the reductions needed to prove such a complexity result inherently use a fixed starting location. A proof that Trail Trap is also PSPACE-complete, as we conjecture in Section~\ref{section:conclusion}, would require some new techniques which we suspect would be of independent interest.

\section{Preliminaries}\label{sec preliminaries} 

We briefly introduce definitions and notation pertaining to degrees, distance, and moves. We adopt the standard graph-theoretic notation used in \cite{Diestel} and refer the reader there for additional background. 

Given a graph $G$, the \emph{degree} of a vertex $v$ in $V(G)$, denoted by $\deg_G(v)$, is the number of edges in $E(G)$ incident to $v$. The \emph{distance} between vertices $u, v \in V(G)$, denoted by $d_G(u,v)$, is the length of the shortest path from $u$ to $v$. By convention, $d_G(u,v) = \infty$ if there is no path from $u$ to $v$ in $G$. We write $\deg(v)$ and $d(u,v)$ if the graph is clear from context. 

In a game of {\sc Trail Trap} on $G$, a \emph{move} $m$ is formally an ordered pair $u \rightarrow v$, where $uv\in E(G)$. Let $u$ be the \emph{tail} of move $m$ and $v$ be the \emph{head}. We let $m^{-1}$ denote the move $v \rightarrow u$. 
We label turns by alternating players, such that \poneposs{} first move is turn $1$, \ptwoposs{} first move is turn $2$, and so on. (There is no formal difference between a move and a turn, but the colloquial usage of turn as a ``time to make a move" is useful for us.) For a natural number $r\leq |E(G)|$, an {\em $r$-partial game} $S$ of {\sc Trail Trap} is an ordered $r$-tuple of moves, such that each edge of $G$ occurs at most once and, for all $i \in \{1, \ldots, r-2\}$, the head of move $i$ is equal to the tail of move $i + 2$. An $r$-partial game is \emph{incomplete} (has a legal continuation) if there exists an edge of $G$ not in $S$ and incident to the head of move $r-1$. 

We also introduce several definitions related to graph connectivity and trails. A graph $G$ is \emph{connected} if there exists a trail in $G$ between any pair of its vertices, and $G$ is \emph{disconnected} if it contains at least two \emph{components}, or maximal connected subgraphs. Given a trail $t = v_1 v_2 \ldots v_{m}$ on $G$, let $\MaxTrail(t) = m-1$ be the length of $t$. Let $\MaxTrail(G)$ be the maximum length of a trail in $G$. We denote the maximum length of a trail in $G$ that starts at a given vertex $v$ by $\MaxTrail(G,v)$.

In {\sc Trail Trap}, each player's trail is restricted to a single component of $G$. On a disconnected graph $G$, there are two ways a game of {\sc Trail Trap} can go: either both players make trails on the same component, or each plays on a separate component. If both players draw trails on the same component $H$, then it suffices to determine a winner for $H$. If both players use different components, however, then the game runs in parallel, with the players unable to interfere with one another's trails. The longest trail in a disconnected graph is important, but does not completely determine who wins. Figure~\ref{fig:disconnected} shows an example graph $G$ where the longest trail is in a component which induces a \pone{}-win graph, but $G$ itself is \ptwo{}-win. In a winning strategy for \pone{} on $G$, by Proposition \ref{lem:p1 loses on degree 1 and 2}, \pone{} must end their first move on the vertex of degree $3$. However, \pone{} can then create a trail of at most four edges, which \ptwo{} can match on the other component. 

\begin{figure}[ht]
    \centering
    \begin{tikzpicture}
    [very thick,every node/.style={circle, draw=black, fill=black, inner sep=2}]

    \node  (0) at (-0.65*3, 0) {};
    \node  (1) at (-0.65*4, 0.65*1) {};
    \node  (2) at (-0.65*5, 0.65*2) {};
    \node  (3) at (-0.65*2, 0.65*1) {};
    \node  (4) at (-0.65*1, 0.65*2) {};
    \node  (5) at (-0.65*3, -0.65*1.25) {};
    \node  (6) at (-0.65*3, -0.65*2.5) {};
    \node  (7) at (0, 0) {};
    \node  (8) at (0.65*1.5, 0) {};
    \node  (9) at (0.65*3, 0) {};
    \node  (10) at (0.65*4.5, 0) {};
    \node  (11) at (-0.65*6, 0.65*3) {};
    \node  (12) at (0, 0.65*3) {};
    \node  (13) at (-0.65*3, -0.65*3.75) {};
    \node  (14) at (0.65*6, 0) {};
	
    \draw (11) to (2);
    \draw (2) to (1);
    \draw (1) to (0);
    \draw (0) to (5);
    \draw (5) to (6);
    \draw (6) to (13);
    \draw (0) to (3);
    \draw (3) to (4);
    \draw (4) to (12);
    \draw (7) to (8);
    \draw (8) to (9);
    \draw (9) to (10);
    \draw (10) to (14);

    \end{tikzpicture}
    \caption{A disconnected \ptwo{}-win graph in which the component containing the longest trail is a \pone{}-win graph.}
    \label{fig:disconnected}
\end{figure}

\begin{proposition}
\label{prop:disconnected}
    A graph $G$ is \pone{}-win if and only if there exists a component $H$ with a trail $t=v_1 v_2 \ldots v_{m+1}$ such that there exists a winning strategy for \pone{} on $H$ with first move $v_1\to v_2$, and for all other components $H'$ of $G$, $\MaxTrail(H')<m$.    
\end{proposition}
\begin{proof}    
    The result is trivial where $G$ is connected, so assume $G$ is disconnected. Suppose first that $G$ is \pone{}-win and that \pone{} has a winning strategy on $G$ beginning with the move $v_1\to v_2$ on some component $H$. It follows that $H$ is \pone{}-win and that there is a winning strategy for \pone{} on $H$ beginning with the move $v_1\to v_2$, otherwise \ptwo{} could win on $G$ by also playing on $H$. Suppose that the longest trail in $H$ beginning with $v_1v_2$ has length $m$. If there exists another component $H'$ where $\MaxTrail(H')\ge m$, then \ptwo{} can play on $H'$ and create a trail at least as long as \pone{}'s longest possible trail. Hence, \ptwo{} would win on $G$. Thus, $\MaxTrail(H')<m$.
    
    Conversely, suppose that there exist: (i) a component $H$ such that $H$ is \pone{}-win, (ii) a winning strategy for $\pone{}$ on $H$ that starts with the move $v_1\to v_2$, and (iii) a trail of length $m$ on $H$ beginning with edge $v_1v_2$, and suppose that no other component $H'$ satisfies $\MaxTrail(H')\ge m$. We can construct a winning strategy for $\pone{}$ on $G$ as follows. First, \pone{} plays the move $v_1\to v_2$. If \ptwo{} plays on a separate component, they get a trail of length less than $m$ while \pone{} can continue along a trail of length $m$. Thus, \ptwo{} is forced to play on $H$. However, as there exists some winning strategy for \pone{} on $H$, beginning with the move $v_1\to v_2$, \pone{} can now win by following that \pone{}-win strategy for {\sc Trail Trap} on $H$.
\end{proof}

To determine the winner on a disconnected graph, we must be able to compare the maximum trail length of components in general and from any given starting vertex. Upper and lower bounds are known for maximum trail length on a graph with a given number of vertices and edges \cite{BOLLOBAS198251}, but despite the invariant's simplicity, little else is known. For a given connected graph $G$, it is well known that there is a trail containing every edge in $G$ if and only if $G$ contains at most two vertices of odd degree. If $G$ contains no vertices of odd degree, then $G$ is Eulerian, and there exists a trail containing every edge and beginning and ending on the same vertex. Some related material exists on determining the number of edges by which a graph is deficient from being Eulerian \cite{BoSuTi77}. Without an independent means of determining maximum trail length, Proposition \ref{prop:disconnected} is insufficient for determining a winner on a disconnected graph. We further show in Section~{\ref{sec NP-hardness}} that determining whether a graph is \ptwo{}-win is NP-hard for a certain class of connected graphs, and thus the problem is NP-hard for disconnected graphs as well.
 
There are some special cases for which the winner of {\sc Trail Trap} can be easily determined. Trivially, if $|E(G)| = \emptyset$, then $G$ is \ptwo{}-win (\pone{} is the first to be unable to move), and if $|E(G)| = 1$, then $G$ is \pone{}-win. We assume henceforth that $|E(G)| \ge 2$. Additionally, a \pone{}-win graph must contain a vertex of degree at least three. 

\begin{proposition}
\label{lem:p1 loses on degree 1 and 2}
    Let $G$ be a \pone{}-win graph with at least two edges. If $u \to v$ is \poneposs{} first move in a winning strategy, then $\deg(v) \geq 3$.
    In particular, the maximum degree of $G$ is at least $3$.
\end{proposition}
\begin{proof}
    Let $G$ be a graph with at least two edges. We denote \poneposs{} first move by $u \to v$ and \ptwoposs{} by $x \to y$.
    Clearly, if $\deg(v) = 1$, then \pone{} cannot make a second move and loses.
    If $\deg(v) = 2$, then $G$ is again a \ptwo{}-win graph, as demonstrated by setting $x = v$ and $y \in N(v) - u$.
\end{proof}

As a result, graphs with maximum degree at most $2$ and at least two edges (including all nontrivial paths and cycles) are \ptwo{}-win. 

Another broad strategy is what Berlekamp, Conway, and Guy \cite{BeCoGu82} call the ``Tweedledee and Tweedledum" strategy, where \ptwo{} mirrors each of \pone{}'s moves. The result also holds on multigraphs. An \emph{automorphism} of a graph $G$ is a bijection $\phi: V(G) \rightarrow V(G)$ that preserves edge relationships, i.e., for $v_i, v_j \in V(G)$, we have $v_iv_j \in E(G)$ if and only if $\phi(v_i)\phi(v_j) \in E(G)$. We call $\phi$ \emph{involutive} if $\phi(\phi(v)) = v$ for all $v \in G$. An edge $v_iv_j$ of $G$ is a \emph{fixed edge} of $\phi$ if the edge is mapped to itself; that is, either $\phi(v_i) = v_j$ and $\phi(v_j) = v_i$, or $\phi(v_i) = v_i$ and $\phi(v_j) = v_j$.

\begin{proposition}
\label{prop:copycat_involution_strategy}
    If a graph $G$ has an involutive automorphism with no fixed edges, then $G$ is \ptwo{}-win.
\end{proposition}
\begin{proof}
    Let $G$ be a graph with an involutive automorphism $\phi: V(G) \rightarrow V(G)$ with no fixed edges. If \pone{} moves $v_{i} \to v_j$ on turn $k$, then \ptwo{} moves $\phi(v_{i}) \to \phi(v_j)$ on turn $k+1$. The move is legal since, if $k \ge 3$, \ptwo{} ended their previous move on $\phi(v_{i})$, edge $v_{i}v_j$ was unused before turn $k$, and $\phi$ fixes no edges. It follows that $\phi(v_{i})\phi(v_j)$ was unused before turn $k+1$. Hence \ptwo{} can always play after \pone{}; as $G$ is finite, \ptwo{} wins. 
\end{proof}

\section{The Game on Common Graph Families}\label{sec solved families}

We consider {\sc Trail Trap} on a number of basic classes of graphs in order to illustrate a variety of winning strategies as well as to point to the difficulty of the problem. Almost all graphs on four or fewer vertices (still with at least two edges) are \ptwo{}-win. There are two exceptions: $K_{1,3}$ and $K_4$ (the former is trivial and the latter is shown in Proposition \ref{prop:K_4}). This trend appears to continue: we show computationally\footnote{Source code for computations throughout the paper: \url{https://github.com/mackenziecarr/Trail-Trap}} that most small connected graphs are \ptwo{}-win (see Table~\ref{tab:small-graphs}). Due to their rarity, we suspect that \pone{}-win graphs may exhibit interesting graph-theoretic properties.

\begin{table}[htb]
\begin{center}
\begin{tabular}{ c c c c}
 $n$ & Number of Graphs & \ptwo{}-win & Percent \ptwo{}-win \\ 
 3 & 2 & 2 & 100\% \\  
 4 & 6 & 4 & 67\% \\  
 5 & 21 & 17 & 81\% \\  
 6 & 112 & 88 & 79\% \\  
 7 & 853 & 734 & 86\% \\  
\end{tabular}
\end{center}
    \caption{The number of \ptwo{}-win graphs out of connected graphs with a small number of vertices.}
    \label{tab:small-graphs}
\end{table}

We use a depth-first recursive algorithm to evaluate if a given graph $G$ is \pone{}-win or \ptwo{}-win. We say a player wins an incomplete $r$-partial game $I$ if they can win the game on $G$ with optimal play, given that the sequence of both players' starting moves constituted the partial game $I$. For any incomplete $r$-partial game $I$, we iterate through all possible moves for the player whose turn it is and assert that they can win $I$ with optimal play if and only if they can win at least one of the resulting $(r+1)$-partial games extending $I$. 
For highly symmetric graphs, such as complete or complete bipartite graphs, we speed up the computations by limiting the choices for the first few moves, without loss of generality.

We continue by considering complete graphs. For small cases, we are able to determine the winner by hand. In the following proposition, the result that $K_4$ is \pone{}-win is due to Arneman, Catanzaro, and Danison \cite{ArCaDa21}, but restated here with a shorter proof. 

\begin{proposition}
\label{prop:K_4}
    The complete graph $K_4$ is \pone{}-win, and $K_5$ is \ptwo{}-win. 
\end{proposition}
\begin{proof}
    Denote the vertices of $K_n$ by $v_1,v_2,\dots,v_n$. On $K_4$, suppose without loss of generality that \pone{} initially makes the move $v_1 \to v_2$.  We consider three cases based on \ptwoposs{} first move: 
    \begin{itemize}
    \item [(i)] If \ptwoposs{} first move also ends at $v_2$, then \pone{} makes their second move on the only unused edge incident to $v_2$, to win.
    \item [(ii)] If \ptwoposs{} first move ends at $v_3$ (respectively, $v_4$), then \pone{} should move $v_2 \to v_4$ (resp. $v_2 \to v_3$). Irrespective of \ptwoposs{} second move, \pone{} moves $v_4 \to v_1$ (resp. $v_3 \to v_1$) to win. 
    \item [(iii)] If \ptwoposs{} first move is $v_3 \to v_1$ (resp., $v_4 \to v_1$) then \pone{} should move $v_2 \to v_4$. \ptwo{} must then move $v_1 \to v_4$ (resp. $v_1 \to v_3$) and \pone{} moves $v_4 \to v_3$ to win.   
    \end{itemize}   

    Therefore, \pone{} has a winning strategy on $K_4$. 
    
    On $K_5$, suppose again without loss of generality that \pone{} initially moves $v_1 \to v_2$. Have \ptwo{} move $v_2 \to v_3$ and then, without loss of generality, \pone{} responds $v_2 \to v_4$. Have \ptwo{} move $v_3 \to v_1$. We consider three cases based on \poneposs{} third move: 
    
    \begin{itemize}
    \item [(i)] If \pone{} now moves $v_4 \to v_1$, then \ptwo{} moves $v_1 \to v_5$ to win.  
    \item [(ii)] If \pone{} now moves $v_4 \to v_3$, then \ptwo{} moves $v_1 \to v_5$. \pone{} must now move $v_3 \to v_5$, and then \ptwo{} moves $v_5 \to v_4$. \pone{} must now move $v_5 \to v_2$, and then \ptwo{} moves $v_4 \to v_1$ to win. 
    \item [(iii)] If \pone{} now moves $v_4 \to v_5$, then \ptwo{} moves $v_1 \to v_4$. Irrespective of \poneposs{} next move, \ptwo{} moves $v_4 \to v_3$ to win. \hfill \qedhere
    \end{itemize}
\end{proof}

For larger complete graphs and other families, analyzing games by hand is prohibitively difficult. We determine computationally that $K_6$ and $K_7$ are \pone{}-win, and $K_8$ and $K_9$ are \ptwo{}-win.

We now turn to complete bipartite graphs, for which we are able to determine a winner on some infinite subfamilies. 

\begin{theorem}
\label{prop:complete_bipartite}
    The complete bipartite graph $K_{p,q}$ is \ptwo{}-win if at least one of $p, q$ is even. If $p=3$ and $q \geq 5$ is odd, then $K_{p,q}$ is \pone{}-win.
\end{theorem}
\begin{proof}
Given $G = K_{p,q}$, let $(L,R)$ be the bipartition of $V(G)$ with $L=\{u_1,\dots,u_p\}$ and $R=\{v_1,\dots,v_q\}$. Consider the involutive automorphism $\phi$ which sends $u_i$ to $u_{p+1-i}$ and $v_j$ to $v_{q+1-j}$. The only way an edge is fixed is if the edge $u_iv_j$ satisfies $i=p+1-i$ and $j=q+1-j$. However, as $i$ and $j$ are integers, this only occurs if $p$ and $q$ are both odd. Thus, if either $p$ or $q$ is even, then $K_{p,q}$ is \ptwo{}-win by Proposition \ref{prop:copycat_involution_strategy}.

For $p=3$ and $q \ge 13$ odd, we show that \pone{} has a winning strategy. Let \poneposs{} first move be $u_1 \to v_1$; after a series of initial moves, whenever \pone{} is moving from $L$ to $R$, they follow the strategy of moving to a vertex $v_j$ in $R$ that has not been previously visited by either player, as long as this is possible. This ensures that even if \ptwo{} immediately uses an edge incident to $v_j$, \pone{} will be able to move back to $L$ via an unused edge.

Suppose first that \ptwo{} initially also moves from $L$ to $R$. This means \ptwo{} is always moving to the partite set ($L$ or $R$) currently containing \poneposs{} token. Before any move by \ptwo{} from $L$ to $R$, every vertex in $R$, except for \pone{}'s current vertex, has an odd number of unused edges. If \ptwo{} moves to a vertex with one unused edge, then they can make no further moves, and \pone{} can win by moving back to $L$ along an unused edge incident to their current vertex. If \ptwo{} moves to a vertex $v_j$ with two unused edges, this must be \pone{}'s current vertex, so then \poneposs{} subsequent move uses the last available edge incident to $v_j$ to win. Thus we may assume that whenever \ptwo{} moves to $R$, they move to a vertex with all three edges previously unused.

According to \pone{}'s strategy of moving to vertices in $R$ with three incident unused edges whenever possible, the two players are taking turns visiting previously unvisited vertices in $R$. Since $|R| = q$ is odd, \ptwo{} will run out of previously unvisited vertices in $R$ to move to right after \pone{} moves to their $\frac{(q+1)}{2}^{\text{th}}$ vertex in $R$. Thus, \ptwo{} will be forced to move to a vertex in $R$ with one or two unused edges, ensuring their defeat as described earlier.

Suppose instead that \ptwo{} initially moves from $R$ to $L$. As long as there are vertices in $R$ with three unused edges, \pone{} can move to such a vertex when going from $L$ to $R$ and subsequently return to $L$ along an unused edge. As such, \pone{} cannot lose until all vertices from $R$ have been visited. We show that \pone{} can force a game position where both \pone{} and \ptwo{} have their tokens on vertices in $L$, it is \ptwoposs{} turn, all vertices in $R$ have zero, one, or three unused edges, and an even number have three unused edges. Given this, we show next that \pone{} can win by forcing \ptwo{} to move to a vertex of $R$ with one unused edge. 

Without loss of generality, \ptwo{} can either start from $v_1$ or from some other vertex of $R$. If \ptwo{} begins at $v_1$, then without loss of generality, they move to $u_2$. \pone{} can then move to $u_3$. At this point, both players are in $L$, it is \ptwoposs{} turn, and $R$ contains one vertex with zero unused edges and $q-1$ vertices with three unused edges.  

If otherwise \ptwo{} begins at another vertex of $R$ (say, $v_2$), then up to symmetry we assume \ptwo{} moves to $u_1$ or $u_2$. \pone{} can then move to $u_3$. If \ptwo{} now moves to $v_1$ (note this is only possible if \ptwoposs{} first move was to $u_2$), they cannot move further as \pone{} has already used $u_1v_1$ and $v_1u_3$, and \pone{} makes another move to win. \ptwo{} cannot legally move to $v_2$, so up to symmetry suppose they move to $v_3$. \pone{} can then move to $v_2$. Suppose that both \ptwo{} and \pone{} now return to any sensible (not immediately losing) vertices in $L$. At this point, it is \ptwoposs{} turn and $R$ contains one vertex with no unused edges, two vertices with one unused edge, and $q-3$ vertices with three unused edges. 

In either case, it is now impossible for \ptwo{} to move to a vertex in $R$ with two unused edges at any future point. Since $q$ is by assumption odd, there are an even number of vertices left in $R$ with three unused edges. It is \ptwo{}'s turn, so ultimately, \ptwo{} will be forced to be first to visit a vertex in $R$ with one unused edge, assuming \pone{} uses the strategy of moving to a vertex in $R$ with three unused edges whenever possible. Once this happens, \ptwo{} cannot make another move and \pone{} can win if they are able to make one final move from $L$ to $R$.

Next, we show that we can specify \poneposs{} midgame strategy to ensure they have a final (winning) move available, i.e., that when \ptwo{} enters a vertex in $R$ along its last unused edge, \pone{} can still make one more move. Suppose \ptwo{} enters a vertex in $R$ along its last unused incident edge on turn $r$ of the game. Necessarily, $r$ is even, and since there are at most $2q+1$ previously used edges, we have $r \le 2q+2$. We also know \pone{}'s previous move was from $R$ to $L$ and they had a choice between two vertices in $L$. Hence as long as there are at least four unused edges incident to these two vertices at the conclusion of turn $r-2$, \pone{} can choose to move either to a vertex in $L$ with three unused edges or to one with only two unused edges and not occupied by \ptwo{}, and thus be able to make one final move to win.

We now show that for odd $q \ge 13$, \pone{} can move in such a way as to guarantee that all three vertices in $L$ have at least four used edges at the conclusion of turn $2q$. This requires at least $12$ moves by \pone{}, so we need $q \ge 12$ such that \pone{} does not run out of unused vertices in $R$. 

Wes specify the following strategy: any move by \pone{} from $L$ to $R$ after turn $7$ but before turn $r$ is to a vertex in $R$ with three unused edges. Such a vertex will still have two unused edges after \ptwo{} makes their next move, so \pone{} will have two choices for where to go next in $L$ and can choose the one they have previously visited fewer times, regardless of what \ptwo{} does. Based on the strategy outlined for the first seven moves of the game, \pone{} begins at $u_1$, then next enters $L$ by moving to $u_3$, then next enters $L$ by moving to either $u_1$ or $u_2$. Since after this point, \pone{} can enter $L$ at either of the two vertices they did not most recently visit, $\pone{}$ can force their first seven visits to $L$ to be either $u_1,u_3,u_2,u_1,u_3,u_2,u_1$ or $u_1,u_3,u_1,u_2,u_3,u_2,u_1$. Counting both the initial edge out of $u_1$ and the final edge into $u_1$ means that $\pone{}$ can use four edges incident to each vertex in $L$ within their first twelve moves, and thus within the first $24< 2q$ total moves. Let \pone{} follow this strategy. 

Thus, for each pair of vertices in $L$, at least four of the first $2q$ moves were incident to the third vertex in $L$, so this pair has at least four unused incident edges remaining after the first $2q$ moves. Since $r-2\le 2q$, any pair has at least four unused incident edges remaining after $r-2$ moves, as desired.

Therefore, $K_{3,q}$ is \pone{}-win when $q$ is odd and $q \ge 13$. We confirm computationally that $K_{3,5},K_{3,7},K_{3,9},$ and $K_{3,11}$ are \pone{}-win, whereas $K_{3,3}$ is \ptwo{}-win.
\end{proof}

If both sides of the bipartition have a large odd number of vertices, the winner remains unknown. Our computations show mixed results: $K_{5, 5}$, $K_{5, 7}$, $K_{5, 11}$, and $K_{7,7}$ are \ptwo{}-win, and $K_{5,9}$ is \pone{}-win. 

In preparation for our next class of graphs, we modify Proposition \ref{prop:copycat_involution_strategy} for $r$-partial games. An $r$-partial game is not equivalent to a full game of {\sc Trail Trap} on the remaining unused portions of the graph, because the players have fixed starting vertices, though not yet edges. To that end, we require the involution on the remaining part of the graph to respect these fixed starting vertices. 

\begin{proposition}
\label{prop:CCI_revisited}
     Given a graph $G$ and an $r$-partial game $S$ on $G$, let $E_S$ be the set of edges of $G$ unused in $S$, let $G_S = (V(G), E_S)$, and let $v_1$ and $v_2$ be the ending vertices of moves $r-1$ and $r$, respectively, of $S$. If there exists an involutive automorphism $\phi: V(G_S) \rightarrow V(G_S)$ of $G_S$ with no fixed edges and $\phi(v_1) = v_2$, then \pone{} wins with optimal play if $r$ is odd and \ptwo{} wins with optimal play if $r$ is even. 
\end{proposition}
\begin{proof}
    The proof is analogous to that of Proposition \ref{prop:copycat_involution_strategy}. Let Player A be the player whose turn it currently is (i.e., if $r$ is even, Player A is \pone, and if $r$ is odd, Player A is \ptwo.). We propose a strategy for Player B. If Player A moves $v_{i}\to v_j$ on turn $r+(2k+1)$, $k \ge 0$, then Player B responds with $\phi(v_{i})\to \phi(v_j)$ on turn $r + (2k + 2)$. For $k = 0$, the move is legal, since by assumption Player A has begun on $v_1$ and player $B$ begins on $v_2 = \phi(v_1)$. Since $\phi$ has no fixed edges, edge $\phi(v_i)\phi(v_j)$ is distinct from edge $v_iv_j$. Suppose that for some $m$, all of Player B's moves $r+(2k+2)$ for $k < m$ are legal, and that Player A has a legal move $v_iv_j$ on turn $r+2m+1$. Then the move $\phi(v_i)\phi(v_j)$ on turn $r+2m+2$ is legal since:
    \begin{enumerate}
        \item [(i)] By assumption, Player A's move $r+2m - 1$ ends on some vertex $v_{i}$, so move $r+2m$ ends on $\phi(v_{i})$.
        \item [(ii)] If the edge $\phi(v_{i})\phi(v_j)$ was used by player B before turn $r+2m+2$, then edge $\phi(\phi(v_i))\phi(\phi(v_j)) = v_iv_j$ was used on the turn prior, and hence before turn $r+2m+1$, for a contradiction. If $\phi(v_{i})\phi(v_j)$ was used previously by Player A, then $v_iv_j$ would have been played by Player B instead of Player A.
    \end{enumerate}
    Since Player B is always able to make a legal move after each legal move of Player A, it follows that this is a winning strategy for player B for the graph $G_S$ and the given starting positions. Therefore, \pone{} can win on $G$ with optimal play if $r$ is odd and \ptwo{} can win on $G$ with optimal play if $r$ is even.
\end{proof}

We next consider graph products. The \emph{Cartesian product} of graphs $G$ and $H$, denoted by $G\,\square \,H$ is the graph with vertex set $V(G\,\square \,H) = V(G) \times V(H)$ where two vertices $(x,y)$ and $(u,v)$ are adjacent in $G\,\square \,H$ if and only if either (i) $x=u$ in $G$ and $yv\in E(H)$; or (ii) $y=v$ in $H$ and $xu\in E(G)$ holds. For a set of edges $S$, let $G-S$ denote the subgraph of $G$ formed by removing the edges in $S$.

Both grid graphs (the Cartesian product of two paths) and prism graphs (the Cartesian product of a cycle and an edge) are solvable in particular cases. For a positive integer $n$, let $P_n$ denote a path graph on $n$ vertices and $n-1$ edges. 
We note that if $m$ and $n$ have the same parity, then $P_m \mathbin{\square} P_n$ has an involutive automorphism with no fixed edges, and is thus \ptwo{}-win by Proposition~\ref{prop:copycat_involution_strategy}.
In what follows, we prove that $P_2 \mathbin{\square} P_n$ is \pone{}-win whenever $n > 3$ is odd. In Section~\ref{section:conclusion}, we conjecture that these are the only \pone{}-win grid graphs.

\begin{proposition}
    \label{prop:grids}
    If $n>3$ is odd, then $P_2 \,\square \, P_n$ is \pone{}-win. 
\end{proposition}
\begin{proof}
    Let $n \ge 3$ be odd, let $k = \frac{n-1}{2}$, and let $G = P_2 \,\square \, P_n$ with $V(G) = \{ u_{-k},u_{-k+1},\dots,u_0, \ldots, u_k, $ $ v_{-k},v_{-k+1},\dots,v_0, \ldots, v_k\}$ and $E(G) = \{u_iv_i : -k \le i \le k \} \cup \{u_i u_{i+1}:-k \le i \le k-1 \} \cup \{ v_i v_{i+1} : -k \le i \le k-1 \}$, as shown in Figure \ref{fig:gridex}. Let \pone{} begin with move $u_0 \to v_0$. We consider each of \ptwoposs{} possible first moves, in turn, assuming up to symmetry that the endpoints have nonnegative index. With one exception (case (viii)), let \poneposs{} second move be $v_0 \to v_{-1}$. We record the general strategies in Table \ref{table:gridgraphcases}. 

\begin{figure}[ht]
    \centering
    \begin{tikzpicture}[very thick,every node/.style={circle, draw=black, fill=black, inner sep=2}]
		\node[label={[label distance=-1mm]$u_{-k}$}] (0) at (-4*1.5, 1.5) {};
		\node [label={[label distance=-2.5mm]$u_{-k+1}$}] (1) at (-3*1.5, 1.5) {};
            \node[draw=none,fill=none] (dots) at (-1.5*1.5,0.5*1.5){$\dots$};
		\node [label={[label distance=-2.5mm]$u_{-k+2}$}] (2) at (-2*1.5, 1.5) {};
		\node[label={[label distance=-0.5mm]$u_{-2}$}]  (3) at (-1.5, 1.5) {};
		\node [label={[label distance=-0.5mm]$u_{-1}$}] (4) at (0, 1.5) {};
		\node [label={[label distance=1mm]$u_0$}]  (5) at (1.5, 1.5) {};
		\node [label={[label distance=1mm]$u_1$}] (6) at (2*1.5, 1.5) {};
		\node [label={[label distance=1mm]$u_2$}] (7) at (3*1.5, 1.5) {};
            \node[draw=none,fill=none] (dots2) at (3.5*1.5,0.5*1.5) {$\dots$};
		\node [label={[label distance=-1mm]$u_{k-2}$}] (8) at (4*1.5, 1.5) {};
		\node[label={[label distance=-1mm]$u_{k-1}$}]  (9) at (5*1.5, 1.5) {};
		\node[label=below:$v_{-k}$]  (10) at (-4*1.5, 0) {};
		\node [label={[label distance=-1.5mm]below:$v_{-k+1}$}] (11) at (-3*1.5, 0) {};
		\node [label={[label distance=-1.5mm]below:$v_{-k+2}$}] (12) at (-2*1.5, 0) {};
		\node[label=below:$v_{-2}$]  (13) at (-1.5, 0) {};
		\node [label=below:$v_{-1}$] (14) at (0, 0) {};
		\node[label={[label distance=1mm]below:$v_0$}]  (15) at (1.5, 0) {};
		\node [label={[label distance=1mm]below:$v_1$}] (16) at (2*1.5, 0) {};
		\node [label={[label distance=1mm]below:$v_2$}] (17) at (3*1.5, 0) {};
		\node [label={[label distance=-0.5mm]below:$v_{k-2}$}] (18) at (4*1.5, 0) {};
		\node [label={[label distance=-0.5mm]below:$v_{k-1}$}]  (19) at (5*1.5, 0) {};
		\node [label={[label distance=0.5mm]$u_{k}$}] (20) at (6*1.5, 1.5) {};
		\node [label={[label distance=1mm]below:$v_k$}] (21) at (6*1.5, 0) {};
        
  \draw (0) -- (10);
		\draw (1) to (11);
		\draw (2) to (12);
  \draw (13) -- (3);
		\draw (4) to (14);
  \draw (5) -- (15);
		\draw (6) to (16);
		\draw (7) to (17);
		\draw (8) to (18);
		\draw (9) to (19);
  \draw (1) -- (0);
  \draw (1) -- (2);
		\draw (3) to (4);
		\draw (4) to (5);
		\draw (5) to (6);
		\draw (6) to (7);
		\draw (8) to (9);
		\draw (10) to (11);
		\draw (11) to (12);
		\draw (13) to (14);
		\draw (14) to (15);
  \draw (15) -- (16);
  \draw (16) -- (17);
  \draw (18) -- (19);
		\draw (9) to (20);
		\draw (19) to (21);
		\draw (21) to (20);
  \draw (2) to (-2.65,1.5);
  \draw (12) to (-2.65,0);
  \draw (3) to (-1.85,1.5);
  \draw (13) to (-1.85,0);
  \draw (7) to (4.85,1.5);
  \draw (17) to (4.85,0);
  \draw (8) to (5.65,1.5);
  \draw  (18) to (5.65,0);
\end{tikzpicture}

\caption{ $G = P_2 \,\square \, P_n$, with vertices labeled as in the proof of Proposition~\ref{prop:grids}. }
\label{fig:gridex}
\end{figure}
    
\begin{table}[ht]
\centering
\begin{tabular}{ |c||c|c| } 
\hline
Case & \ptwoposs{} First Move & Strategy Overview \\ 
\hline \hline
\multirow{1}{1cm}{(i)} & $v_{1} \to v_{0}$  & \multirow{1}{8.5cm}{\pone{} moves $v_0 \to v_{-1}$ and \ptwo{} is blocked} \\ 
\hline
\multirow{1}{1cm}{(ii)} & $u_0 \to u_{1}$  & \multirow{1}{8.5cm}{Involution $\phi_1: u_{-x} \mapsto v_{x} \text{ and } v_{-x} \mapsto u_{x}$ using Prop. \ref{prop:CCI_revisited}}\\  
\hline
\multirow{1}{1cm}{(iii)} & $v_0 \to v_{1}$  & \multirow{1}{8.5cm}{Involution $\phi_2: u_{-x} \mapsto u_{x} \text{ and } v_{-x} \mapsto v_{x}$ using Prop. \ref{prop:CCI_revisited}}\\ 
\hline
\multirow{1}{1cm}{(iv)} & $u_{1} \to u_0$  & \multirow{1}{8.5cm}{\pone{} moves $v_{x} \to v_{x-1}$ on each turn} \\
\hline
\multirow{1}{1cm}{(v)} & $u_{1} \to v_{1}$  & \multirow{1}{8.5cm}{Involution $\phi_2: u_{-x} \mapsto u_{x} \text{ and } v_{-x} \mapsto v_{x}$ using Prop. \ref{prop:CCI_revisited}} \\
\hline
\multirow{9}{1cm}{(vi)} & $v_{1} \to u_{1}$  & \multirow{9}{8.5cm}{Zigzag strategy (see Figure \ref{fig:zigzag})} \\
& $v_{1} \to v_{2}$  & \\
& $u_{2} \to u_{1}$  & \\
& $u_{2} \to v_{2}$  & \\
& $v_{2} \to v_{1}$  & \\
& $u_{3} \to u_{2}$  & \\
& $v_{3} \to v_{2}$  & \\
& $v_{3} \to u_{3}$  & \\
& $u_{4} \to u_{3}$  & \\
\hline
\multirow{1}{1cm}{(vii)} & $v_{2} \to u_{2}$  & \multirow{1}{8.5cm}{Zigzag strategy or Loop strategy (see Figure \ref{fig:gridsvi})} \\
\hline
\multirow{4}{1cm}{(viii)} & $u_{1} \to u_{2}$ & \multirow{4}{8.5cm}{Blocking strategy (see Figure \ref{fig:gridsvii})}\\
& $u_{2} \to u_{3}$ & \\
& $z \to u_{x}$, $x \ge 4$ & \\
& $z \to v_{x}$, $x \ge 3$ & \\
\hline
\end{tabular}
\caption{Cases for \poneposs{} strategy on the graph $P_2 \,\square \, P_n$, depending on \ptwoposs{} first move, with $z$ designating an arbitrary vertex.}
\label{table:gridgraphcases}
\end{table}
\begin{itemize}
\item [(i)] Let \ptwo{} begin with move $v_{1} \to v_{0}$, which \pone{} follows with $v_0 \to v_{-1}$. \pone{} wins the game as \ptwo{} is unable to make a second move. 

\item [(ii)] Let \ptwo{} begin with move $u_0 \to u_{1}$, which \pone{} follows with $v_0 \to v_{-1}$. Consider the involutive automorphism on $G$ sending $u_{-x}$ to $v_{x}$ and $v_{-x}$ to $u_{x}$, for $-k \le x \le k$ (imagining $G$ as a ladder, this involution represents a $180^\circ$ rotation of the graph about the edge $u_0v_0$). This automorphism meets the conditions of Proposition \ref{prop:CCI_revisited}, so \pone{} wins on $G$. 

\item [(iii)] Let \ptwo{} begin with move $v_{0} \to v_{1}$, which \pone{} follows with $v_0 \to v_{-1}$. Consider the involutive automorphism on $G$ sending $u_{-x}$ to $u_{x}$ and $v_{-x}$ to $v_{x}$, for $-k \le x \le k$ (imagining $G$ as a ladder, this involution represents a vertical flip of the graph about the edge $u_0v_0$). This automorphism meets the conditions of Proposition \ref{prop:CCI_revisited}, so \pone{} wins on $G$. 

\item [(iv)] Let \ptwo{} begin with move $u_{1} \to u_0$, which \pone{} follows with $v_0 \to v_{-1}$. On their next turn, \ptwo{} must move $u_0 \to u_{-1}$. On each subsequent turn, \pone{} should move $v_{x} \to v_{x-1}$. If \ptwo{} ever moves $u_{x} \to v_{x}$, then since \pone{} has already played on the edge $v_{x} v_{x-1}$, \ptwo{} will be unable to make another move. Thus \pone{} wins the game after move $v_{x-1}\to u_{x-1}$, which is still available. Otherwise, \ptwo{} ultimately moves $u_{-k+1} \to u_{-k}$, and \pone{} moves $v_{-k} \to u_{-k}$ to win the game. 

\item [(v)] Let \ptwo{} begin with move $u_{1} \to v_{1}$, which \pone{} follows with $v_0 \to v_{-1}$. Let $E_S$ be the set of unused edges in $G$ after these first three moves, and let $G_S = (V(G), E_S)$. Consider the involutive automorphism on $G_S-\{u_{-1}v_{-1},v_0v_1\}$ sending $u_{-x}$ to $u_{x}$ and $v_{-x}$ to $v_{x}$, for $-k \le x \le k$ (imagining $G$ as a ladder, this involution represents a vertical flip of the graph about the edge $u_0v_0$). This automorphism meets the conditions of Proposition \ref{prop:CCI_revisited}, so \pone{} wins on $G_S-\{u_{-1}v_{-1},v_0v_1\}$. Moreover, if \ptwo{} moves on edge $v_1 v_0$ on their second turn, then they immediately lose since \pone{} is able to make a third move but \ptwo{} is not. If \ptwo{} attempts to use edge $u_{-1}v_{-1}$, then at some earlier point they must have made moves $u_1 \to u_0$ and $u_0 \to u_{-1}$. However, by the involutive automorphism, \pone{} would have previously made the move $u_{-1} \to u_0$, invalidating the second move and securing an immediate win for \pone{}. Hence we can assume \ptwo{} uses neither of these edges, and the given strategy on $G_S-\{u_{-1}v_{-1},v_0v_1\}$ is also a winning strategy for \pone{} on $G$. 

\item [(vi)] Let \ptwo{} begin with any of the following moves: $v_{1} \to u_{1}, v_{1} \to v_{2}, u_{2} \to u_{1}, u_{2} \to v_{2}, v_{2} \to v_{1}, u_{3} \to u_{2}, v_{3} \to v_{2}, v_{3} \to u_{3}, u_{4} \to u_{3}$. From here, \pone{} will make the following sequence of moves, which we call the \emph{zigzag strategy} (see Figure~\ref{fig:zigzag} for an illustration). Beginning on the next turn with $\ell=0$, \pone{} will move $v_{-2\ell} \to v_{-2\ell -1}$, then $v_{-2\ell -1} \to u_{-2\ell-1}$, then $u_{-2\ell-1} \to u_{-2\ell-2}$, then $u_{-2\ell -2} \to v_{-2\ell -2}$. On the following turn, increase $\ell$ by 1 and repeat the sequence of moves, until the edge $u_{-k}v_{-k}$ is played (in some direction). \pone{} can make as a final move whichever of $v_{-k} v_{-k+1}$ or $u_{-k} u_{-k+1}$ is incident to their current vertex. This creates a trail beginning with $u_0 v_0$ that has length $n+1$. 

Meanwhile, \ptwo{} can only play on edges with the tail among $u_1,u_2,\dots,u_{k}, v_1,v_2,\dots,v_{k}$, as well as possible final move $u_0 \to u_{-1}$ (\ptwo{} cannot go further into the negative-indexed vertices because \pone{} will have already used edges $u_{-1}v_{-1}$ and $u_{-1}u_{-2}$). We note also that \ptwo{} can use at most one of $u_1u_0$ and $v_1v_0$ on their trail. The vertices $u_1,u_2,\dots,u_{k}, v_1,v_2,\dots,v_{k}$ can be the tail of exactly one move (since each has degree at most $3$), save for the starting vertex, which may be the tail of two moves. In order for \ptwo{} to build a trail with $n+1=2k+2$ moves, they must secure a trail of $n$ edges whose tails are among vertices $u_1,u_2,\dots,u_{k}, v_1,v_2,\dots,v_{k}$, then end their trail with move $u_0 \to u_{-1}$. However, as the graph is bipartite and such a trail has even length, this necessitates that \ptwo{} started on a vertex in the same partite class as $u_{-1}$, namely $u_{2\ell+1}$ or $v_{2\ell}$ for some integer $\ell$. 

Thus, the only remaining options for \ptwoposs{} first move are $v_2\to v_1$ and $u_3\to u_2$. \ptwoposs{} move before $u_0\to u_{-1}$ must be $u_1\to u_0$, so we require a trail among vertices $u_1,u_2,\dots,u_k, v_1,v_2,\dots,v_k$ that contains the starting vertex twice and the others once, ending with $u_{1}$.
With starting move $u_{x} \to u_{x-1}$ (respectively, $v_{x} \to v_{x-1}$), the parallel edge $v_{x}v_{x-1}$ (respectively, $u_{x}u_{x-1}$) must be used in the reverse direction $v_{x-1} \to v_{x}$, but then the this series of moves cannot end at $u_{1}$. Therefore, \ptwo{} cannot build a trail of length $n+1$, and \pone{} wins after a game of at most $2n+1$ turns.

        \begin{figure}[ht]
        \centering
            \begin{tikzpicture} [very thick,every node/.style={circle, draw=black, fill=black, inner sep=2}]
		\node [label={[label distance=-3mm]$u_{-2\ell-2}$}] (0) at (-4, 1.5) {};
		\node[label={[label distance=-3mm]$u_{-2\ell-1}$}] (1) at (-2, 1.5) {};
		\node[label={[label distance=-1.5mm]$u_{-2\ell}$}] (2) at (0, 1.5) {};
		\node[label={[label distance=-3mm]below:$v_{-2\ell-2}$}] (4) at (-4, 0) {};
		\node[label={[label distance=-3mm]below:$v_{-2\ell-1}$}] (5) at (-2, 0) {};
		\node[label={[label distance=-1.5mm]below:$v_{-2\ell}$}] (6) at (0, 0) {};

   \draw[-{Latex[length=3mm]},color=red] (0) -- node[draw=none,fill=none,left]{\pone{}}(4);
  \draw[-{Latex[length=3mm]},color=red] (5) -- node[draw=none,fill=none,right]{\pone{}}(1);
  \draw[-{Latex[length=3mm]},color=red] (2) -- node[draw=none,fill=none,right]{\pone{}}(6);
		\draw (4) to (5);
  \draw[-{Latex[length=3mm]},color=red] (1) -- node[draw=none,fill=none,above]{\pone{}}(0);
		\draw (2) to (0.5,1.5);
  \draw[-{Latex[length=3mm]},color=red] (6) -- node[draw=none,fill=none,above]{\pone{}}(5);
		\draw (1) to (2);
		\draw (6) to (0.5,0);
            \draw (0) to (-4.5,1.5);
            \draw (4) to (-4.5,0);
 
\end{tikzpicture}

\caption{A sequence of five moves in the zigzag strategy for \pone{} in the proof of Proposition~\ref{prop:grids}.}
\label{fig:zigzag}
\end{figure}

\item [(vii)] Let \ptwo{} begin with move $v_{2}\to u_{2}$, which \pone{} follows with $v_0 \to v_{-1}$. On their next turn, \ptwo{} can play either $u_2\to u_1$, or, if $k\ge 3$, $u_2\to u_3$. If \ptwo{} continues with $u_{2}\to u_{1}$, \pone{} will play $v_{-1}\to u_{-1}$.  Here, \ptwo{} is effectively limited to vertices of index at least $1$. Each of these $2k$ vertices can be the tail of exactly one move (since each has degree at most $3$), save for $v_2$, which may be the tail of two moves. Hence \ptwo{} is able to secure a trail of at most $2k + 1 = n$ edges. At this point, \pone{} will be able to guarantee a trail of length $n+1$ via the zigzag strategy, so \pone{} wins.

If instead \ptwoposs{} second move is $u_{2}\to u_{3}$ (with $k \ge 3$, accordingly), \pone{} will respond on each turn using a strategy that we will call the \emph{loop strategy}. \ptwo{} will continue to make moves of the form $u_{x}\to u_{x+1}$ before finally making a move of the form $u_{x}\to v_{x}$, either by choice or out of necessity. Suppose that this first such move is $u_{a}\to v_{a}$ which would thus be \ptwoposs{} $a^{th}$ move, with $a\ge 3$. Meanwhile, \pone{} is to continue to make moves of the form $v_{-x}\to v_{-x-1}$ so that after $a$ moves, \pone{} is at $v_{-a+1}$. See Figure~\ref{subfig:1} for an illustration of this case. On \poneposs{} $(a+1)^{th}$ turn, they will play $v_{-a+1}\to v_{-a}$.

If \ptwo{} now plays $v_{a}\to v_{a+1}$, \ptwo{} cannot use any further edges with tails less than $a+1$, and thus cannot make more than $n-2a-1$ additional moves, so is ultimately limited to a trail of length at most $(a+1)+(n-2a-1)=n-a\le n-3$. In this case, \pone{} could continue with moves of the form $v_{-x}\to v_{-x-1}$ until reaching $v_{-k}$, then play $v_{-k}\to u_{-k}$, followed by moves of the form $u_{-x}\to u_{-x+1}$ until reaching $u_0$ for a total of $2(k+1) = n+1$ moves. 

Therefore, on \ptwoposs{} $(a+1)^{th}$ turn, suppose that they make the move $v_{a}\to v_{a-1}$. \pone{} follows with $v_{-a}\to u_{-a}$. Now suppose that \ptwo{} keeps making moves of the form $v_{x}\to v_{x-1}$ at least until $v_{2}\to v_{1}$; otherwise they are moving to a vertex with no other unused edges and \pone{} will be able to win just by making an additional move. Meanwhile \pone{} proceeds to make moves of the form $u_{-x}\to u_{-x+1}$. After each player has made $2a-1$ moves, $\pone{}$ is at $u_{-3}$ and \ptwo{} is at $v_{1}$. \pone{} moves $u_{-3}\to u_{-2}$ and \ptwo{} either moves $v_{1}\to u_{1}$ or loses on their next turn. See Figure~\ref{subfig:2} for an illustration of this case. \pone{} now moves $u_{-2}\to u_{-1}$ and regardless of \ptwoposs{} response, \pone{} wins with the move $u_{-1}\to u_0$.

\begin{figure}[ht]
    \centering
    \begin{subfigure}{0.95\textwidth}
    \centering
    \begin{tikzpicture}[very thick,every node/.style={circle, draw=black, fill=black, inner sep=2}]
		\node[label=$u_{a}$] (0) at (-1*-4*1.5, 1.5) {};
		\node  (1) at (-1*-3*1.5, 1.5) {};
            \node[draw=none,fill=none] (dots) at (-1*-2.5*1.5,0.5*1.5){$\dots$};
		\node  (2) at (-1*-2*1.5, 1.5) {};
		\node[label=$u_{2}$]  (3) at (-1*-1.5, 1.5) {};
		\node  (4) at (-1*0, 1.5) {};
		\node [label={[label distance=0.5mm]$u_{0}$}]  (5) at (-1*1.5, 1.5) {};
		\node  (6) at (-1*2*1.5, 1.5) {};
		\node  (7) at (-1*3*1.5, 1.5) {};
            \node[draw=none,fill=none] (dots2) at (-1*3.5*1.5,0.5*1.5) {$\dots$};
		\node  (8) at (-1*4*1.5, 1.5) {};
		\node[label={[label distance=-2.5mm]$u_{-a+1}$}]  (9) at (-1*5*1.5, 1.5) {};
		\node[label=below:$v_{a}$]  (10) at (-1*-4*1.5, 0) {};
		\node  (11) at (-1*-3*1.5, 0) {};
		\node  (12) at (-1*-2*1.5, 0) {};
		\node[label=below:$v_{2}$]  (13) at (-1*-1.5, 0) {};
		\node  (14) at (-1*0, 0) {};
		\node[label={[label distance=0mm]below:$v_0$}]  (15) at (-1*1.5, 0) {};
		\node  (16) at (-1*2*1.5, 0) {};
		\node  (17) at (-1*3*1.5, 0) {};
		\node  (18) at (-1*4*1.5, 0) {};
		\node [label={[label distance=-2.5mm]below:$v_{-a+1}$}]  (19) at (-1*5*1.5, 0) {};
		\node  (20) at (-1*6*1.5, 1.5) {};
		\node  (21) at (-1*6*1.5, 0) {};
        
  \draw[-{Latex[length=3mm]},color=cyan!70!black] (0) -- node[draw=none,fill=none,left]{\ptwo{}}(10);
		\draw (1) to (11);
		\draw (2) to (12);
  \draw[-{Latex[length=3mm]},color=cyan!70!black] (13) -- node[draw=none,fill=none,right]{\ptwo{}}(3);
		\draw (4) to (14);
  \draw[-{Latex[length=3mm]},color=red] (5) -- node[draw=none,fill=none,left]{\pone{}}(15);
		\draw (6) to (16);
		\draw (7) to (17);
		\draw (8) to (18);
		\draw (9) to (19);
  \draw[-{Latex[length=3mm]},color=cyan!70!black] (1) -- node[draw=none,fill=none,above]{}(0);
  \draw[-{Latex[length=3mm]},color=cyan!70!black] (3) -- node[draw=none,fill=none,above]{}(2);
		\draw (3) to (4);
		\draw (4) to (5);
		\draw (5) to (6);
		\draw (6) to (7);
		\draw (8) to (9);
		\draw (10) to (11);
		\draw (12) to (13);
		\draw (13) to (14);
		\draw (14) to (15);
  \draw[-{Latex[length=3mm]},color=red] (15) -- node[draw=none,fill=none,above]{}(16);
  \draw[-{Latex[length=3mm]},color=red] (16) -- node[draw=none,fill=none,above]{}(17);
  \draw[-{Latex[length=3mm]},color=red] (18) -- node[draw=none,fill=none,above]{\pone{}}(19);
		\draw (9) to (20);
		\draw (19) to (21);
		\draw (21) to (20);
  \draw (0) to (-1*-6.5,1.5);
  \draw (10) to (-1*-6.5,0);
  \draw [color=cyan!70!black] (1) to (-1*-4.15,1.5);
  \draw (11) to (-1*-4.15,0);
  \draw [color=cyan!70!black] (2) to (-1*-3.35,1.5);
  \draw (12) to (-1*-3.35,0);
  \draw (7) to (-1*4.85,1.5);
  \draw [color=red] (17) to (-1*4.85,0);
  \draw (8) to (-1*5.65,1.5);
  \draw [color=red] (18) to (-1*5.65,0);
  \draw (20) to (-1*9.35,1.5);
  \draw (21) to (-1*9.35,0);
\end{tikzpicture}
\vspace{-0.75cm}
\caption{Trails built in the first $2a$ turns, with \pone{} starting at $u_0$ and \ptwo{} at $v_{2}$.}
\label{subfig:1}
\end{subfigure}

\begin{subfigure}{0.95\textwidth}
    \centering
    \begin{tikzpicture}[very thick,every node/.style={circle, draw=black, fill=black, inner sep=2}]
		\node[label=$u_{a}$] (0) at (-1*-4*1.5, 1.5) {};
		\node  (1) at (-1*-3*1.5, 1.5) {};
            \node[draw=none,fill=none] (dots) at (-1*-2.5*1.5,0.5*1.5){$\dots$};
		\node  (2) at (-1*-2*1.5, 1.5) {};
		\node[label=$u_{2}$]  (3) at (-1*-1.5, 1.5) {};
		\node  (4) at (-1*0, 1.5) {};
		\node [label={[label distance=0.5mm]$u_{0}$}]  (5) at (-1*1.5, 1.5) {};
		\node  (6) at (-1*2*1.5, 1.5) {};
		\node  (7) at (-1*3*1.5, 1.5) {};
            \node[draw=none,fill=none] (dots2) at (-1*3.5*1.5,0.5*1.5) {$\dots$};
		\node  (8) at (-1*4*1.5, 1.5) {};
		\node[label={[label distance=-2.5mm]$u_{-a+1}$}]  (9) at (-1*5*1.5, 1.5) {};
		\node[label=below:$v_{a}$]  (10) at (-1*-4*1.5, 0) {};
		\node  (11) at (-1*-3*1.5, 0) {};
		\node  (12) at (-1*-2*1.5, 0) {};
		\node[label=below:$v_{2}$]  (13) at (-1*-1.5, 0) {};
		\node  (14) at (-1*0, 0) {};
		\node[label={[label distance=0mm]below:$v_0$}]  (15) at (-1*1.5, 0) {};
		\node  (16) at (-1*2*1.5, 0) {};
		\node  (17) at (-1*3*1.5, 0) {};
		\node  (18) at (-1*4*1.5, 0) {};
		\node [label={[label distance=-2.5mm]below:$v_{-a+1}$}]  (19) at (-1*5*1.5, 0) {};
		\node  (20) at (-1*6*1.5, 1.5) {};
		\node  (21) at (-1*6*1.5, 0) {};
  \draw[-{Latex[length=3mm]},color=cyan!70!black] (0) -- node[draw=none,fill=none,left]{}(10);
		\draw (1) to (11);
		\draw (2) to (12);
  \draw[-{Latex[length=3mm]},color=cyan!70!black] (13) -- node[draw=none,fill=none,right]{\ptwo{}}(3);
  \draw[-{Latex[length=3mm]},color=cyan!70!black] (14) -- node[draw=none,fill=none,right]{\ptwo{}}(4);
  \draw[-{Latex[length=3mm]},color=red] (5) -- node[draw=none,fill=none,left]{\pone{}}(15);
		\draw (6) to (16);
		\draw (7) to (17);
		\draw (8) to (18);
		\draw (9) to (19);
  \draw[-{Latex[length=3mm]},color=cyan!70!black] (1) -- node[draw=none,fill=none,above]{}(0);
  \draw[-{Latex[length=3mm]},color=cyan!70!black] (3) -- node[draw=none,fill=none,above]{}(2);
		\draw (3) to (4);
		\draw (4) to (5);
		\draw (5) to (6);
  \draw[-{Latex[length=3mm]},color=red] (7) -- node[draw=none,fill=none,above]{\pone{}}(6);
  \draw[-{Latex[length=3mm]},color=red] (9) -- node[draw=none,fill=none,above]{}(8);
  \draw[-{Latex[length=3mm]},color=cyan!70!black] (10) -- node[draw=none,fill=none,above]{}(11);
  \draw[-{Latex[length=3mm]},color=cyan!70!black] (12) -- node[draw=none,fill=none,above]{}(13);
  \draw[-{Latex[length=3mm]},color=cyan!70!black] (13) -- node[draw=none,fill=none,above]{}(14);
		\draw (14) to (15);
  \draw[-{Latex[length=3mm]},color=red] (15) -- node[draw=none,fill=none,above]{}(16);
  \draw[-{Latex[length=3mm]},color=red] (16) -- node[draw=none,fill=none,above]{}(17);
  \draw[-{Latex[length=3mm]},color=red] (18) -- node[draw=none,fill=none,above]{}(19);
  \draw[-{Latex[length=3mm]},color=red] (20) -- node[draw=none,fill=none,above]{}(9);
  \draw[-{Latex[length=3mm]},color=red] (19) -- node[draw=none,fill=none,above]{}(21);
  \draw[-{Latex[length=3mm]},color=red] (21) -- node[draw=none,fill=none,right]{}(20);
  \draw (0) to (-1*-6.5,1.5);
  \draw (10) to (-1*-6.5,0);
  \draw [color=cyan!70!black] (1) to (-1*-4.15,1.5);
  \draw [color=cyan!70!black] (11) to (-1*-4.15,0);
  \draw [color=cyan!70!black] (2) to (-1*-3.35,1.5);
  \draw [color=cyan!70!black] (12) to (-1*-3.35,0);
  \draw [color=red] (7) to (-1*4.85,1.5);
  \draw [color=red] (17) to (-1*4.85,0);
  \draw [color=red] (8) to (-1*5.65,1.5);
  \draw [color=red] (18) to (-1*5.65,0);
  \draw (20) to (-1*9.35,1.5);
  \draw (21) to (-1*9.35,0);
\end{tikzpicture}
\vspace{-0.75cm}

\caption{Trails built in the first $4a+1$ turns, with \pone{} starting at $u_0$ and \ptwo{} at $v_{2}$.}
\label{subfig:2}
\end{subfigure}

    \caption{Loop strategy for \pone{} in case (vii) of the proof of Proposition \ref{prop:grids}. Labels indicate the first and most recent moves by each player.}
    \label{fig:gridsvi}
\end{figure}

    \item [(viii)] Let \ptwo{} begin with $u_{1} \to u_{2}, u_{2} \to u_{3}$ or any move whose head has index at least $u_{4}$ or $v_{3}$. Then \pone{} will play the following sequence of moves shown in Figure \ref{fig:gridsvii}, which we call the \emph{blocking strategy}. First \pone{} will play the move series $v_0\to v_{1}, v_{1}\to u_{1}, u_{1}\to u_0$. Note that \ptwo{} is too far away from these edges to claim them before \pone{} can. Specifically, \ptwo{} will not be able to play on edge $u_{1}u_0$ within their first three moves. At this point, \pone{} will continue by playing the following sequence of moves beginning with $\ell=0$. \pone{} will move $u_{-2\ell}\to u_{-2\ell-1}$, then $u_{-2\ell-1}\to v_{-2\ell-1}$, then $v_{-2\ell-1}\to v_{-2\ell-2}$, then $v_{-2\ell-2}\to u_{-2\ell-2}$. On the following turn, increase $\ell$ by 1 and repeat the sequence of moves until edge $u_{-k}v_{-k}$ is played. This guarantees \pone{} a trail of length at least $n+3$. 

    \begin{figure}[ht]
    \centering
    \begin{tikzpicture}[very thick,every node/.style={circle, draw=black, fill=black, inner sep=2}]
		
		\node[label=$u_{2}$]  (3) at (-1*-1.5, 1.5) {};
		\node  (4) at (-1*0, 1.5) {};
		\node [label={[label distance=0.0mm]$u_{0}$}]  (5) at (-1*1.5, 1.5) {};
		\node  (6) at (-1*2*1.5, 1.5) {};
		\node  (7) at (-1*3*1.5, 1.5) {};
            \node[draw=none,fill=none] (dots2) at (-1*3.5*1.5,0.5*1.5) {$\dots$};
		\node [label={[label distance=-1.5mm]$u_{-2\ell}$}]  (8) at (-1*4*1.5, 1.5) {};
		\node  (9) at (-1*5*1.5, 1.5) {};
		
		\node[label=below:$v_{2}$]  (13) at (-1*-1.5, 0) {};
		\node  (14) at (-1*0, 0) {};
		\node[label={[label distance=0mm]below:$v_0$}]  (15) at (-1*1.5, 0) {};
		\node  (16) at (-1*2*1.5, 0) {};
		\node  (17) at (-1*3*1.5, 0) {};
		\node [label={[label distance=-1.5mm]below:$v_{-2\ell}$}] (18) at (-1*4*1.5, 0) {};
		\node   (19) at (-1*5*1.5, 0) {};
		\node [label={[label distance=-3.25mm]$u_{-2\ell-2}$}]  (20) at (-1*6*1.5, 1.5) {};
		\node [label={[label distance=-3.25mm]below:$v_{-2\ell-2}$}] (21) at (-1*6*1.5, 0) {};

		\draw (3) to (13);
        \draw[-{Latex[length=3mm]},color=red] (14) -- node[draw=none,fill=none,above]{}(4);
  \draw[-{Latex[length=3mm]},color=red] (5) -- node[draw=none,fill=none,left]{\pone{}}(15);
        \draw[-{Latex[length=3mm]},color=red] (6) -- node[draw=none,fill=none,above]{}(16);
        \draw[-{Latex[length=3mm]},color=red] (17) -- node[draw=none,fill=none,above]{}(7);
        \draw[-{Latex[length=3mm]},color=red] (18) -- node[draw=none,fill=none,above]{}(8);
        \draw[-{Latex[length=3mm]},color=red] (9) -- node[draw=none,fill=none,above]{}(19);
 
		\draw (3) to (4);
        \draw[-{Latex[length=3mm]},color=red] (4) -- node[draw=none,fill=none,above]{}(5);
        \draw[-{Latex[length=3mm]},color=red] (5) -- node[draw=none,fill=none,above]{}(6);
		\draw (6) to (7);
        \draw[-{Latex[length=3mm]},color=red] (8) -- node[draw=none,fill=none,above]{}(9);
		
		\draw (13) to (14);
        \draw[-{Latex[length=3mm]},color=red] (15) -- node[draw=none,fill=none,above]{}(14);
		\draw (15) to (16);
  \draw[-{Latex[length=3mm]},color=red] (16) -- node[draw=none,fill=none,above]{}(17);
		\draw (18) to (19);
		\draw (9) to (20);
        \draw[-{Latex[length=3mm]},color=red] (19) -- node[draw=none,fill=none,above]{}(21);
        \draw[-{Latex[length=3mm]},color=red] (21) -- node[draw=none,fill=none,right]{\pone{}}(20);
 
  \draw  (3) to (-1*-1.85,1.5);
  \draw (13) to (-1*-1.85,0);
  \draw [color=red] (7) to (-1*4.85,1.5);
  \draw (17) to (-1*4.85,0);
  \draw (8) to (-1*5.65,1.5);
  \draw [color=red] (18) to (-1*5.65,0);
  \draw  (20) to (-1*9.35,1.5);
  \draw (21) to (-1*9.35,0);
\end{tikzpicture}

    \caption{Blocking strategy for \pone{} in case (viii) of the proof of Proposition \ref{prop:grids}. Labels indicate the first and most recent moves by P1. }
    \label{fig:gridsvii}
\end{figure}
    
    Meanwhile \ptwo{} can only play on edges with both endpoints having index at least $1$, which as in cases (vi) and (vii) yields a trail of at most $2k + 1 = n$ edges, for a loss.  \hfill \qedhere
    \end{itemize}
\end{proof}

Let $C_n$ denote the cycle graph on $n$ vertices.

\begin{proposition}
    \label{prop:circularladders}
    The $n$-gonal prism graph, $C_n \,\square \, K_2$, is \ptwo{}-win for all $n$.  
\end{proposition}
\begin{proof}
If $n$ is even, then \ptwo{} has a winning strategy by Proposition~\ref{prop:copycat_involution_strategy}. Now, suppose that $n$ is odd, with $n=2k+1$ for some integer $k$. Let $U=\{u_{-k},$ $u_{-k+1},$ $\dots,u_k\}$ and $V = \{v_{-k},$ $v_{-k+1},$ $\dots,v_k\}$ each induce $C_n$ with vertices in the specified order. Then $G:=C_n \,\square \, K_2$ is formed by taking the disjoint union $U\cup V$ of the two cycles along with the additional edges $u_iv_i$ for $-k \le i \le k$, as shown in Figure \ref{fig:circularladder}. 

\begin{figure}[ht]
\centering
    \begin{tikzpicture}[very thick,every node/.style={circle, draw=black, fill=black, inner sep=2}]
		\node[label={[label distance=-1mm]left:$u_{-k}$}] (0) at (-4*1.5, 1.5) {};
		\node [label={[label distance=-2.5mm]$u_{-k+1}$}] (1) at (-3*1.5, 1.5) {};
            \node[draw=none,fill=none] (dots) at (-1.5*1.5,0.5*1.5){$\dots$};
		\node [label={[label distance=-2.5mm]$u_{-k+2}$}] (2) at (-2*1.5, 1.5) {};
		\node[label={[label distance=-0.5mm]$u_{-1}$}]  (3) at (-1.5, 1.5) {};
		\node [label={[label distance=0.5mm]$u_{0}$}] (4) at (0, 1.5) {};
		\node [label={[label distance=1mm]$u_1$}]  (5) at (1.5, 1.5) {};
		\node [label={[label distance=1mm]$u_2$}] (6) at (2*1.5, 1.5) {};
            \node[draw=none,fill=none] (dots2) at (2.5*1.5,0.5*1.5) {$\dots$};
		\node [label={[label distance=-1mm]$u_{k-2}$}] (8) at (3*1.5, 1.5) {};
		\node[label={[label distance=-1mm]$u_{k-1}$}]  (9) at (4*1.5, 1.5) {};
		\node[label=left:$v_{-k}$]  (10) at (-4*1.5, 0) {};
		\node [label={[label distance=-1.5mm]below:$v_{-k+1}$}] (11) at (-3*1.5, 0) {};
		\node [label={[label distance=-1.5mm]below:$v_{-k+2}$}] (12) at (-2*1.5, 0) {};
		\node[label=below:$v_{-1}$]  (13) at (-1.5, 0) {};
		\node [label={[label distance=1.5mm]below:$v_{0}$}] (14) at (0, 0) {};
		\node[label={[label distance=1mm]below:$v_1$}]  (15) at (1.5, 0) {};
		\node [label={[label distance=1mm]below:$v_2$}] (16) at (2*1.5, 0) {};
		\node [label={[label distance=-0.5mm]below:$v_{k-2}$}] (18) at (3*1.5, 0) {};
		\node [label={[label distance=-0.5mm]below:$v_{k-1}$}]  (19) at (4*1.5, 0) {};
		\node [label={[label distance=0.5mm]right:$u_{k}$}] (20) at (5*1.5, 1.5) {};
		\node [label={[label distance=1mm]right:$v_k$}] (21) at (5*1.5, 0) {};
        
  \draw (0) -- (10);
		\draw (1) to (11);
		\draw (2) to (12);
  \draw (13) -- (3);
		\draw (4) to (14);
  \draw (5) -- (15);
		\draw (6) to (16);
		\draw (8) to (18);
		\draw (9) to (19);
  \draw (1) -- (0);
  \draw (1) -- (2);
		\draw (3) to (4);
		\draw (4) to (5);
		\draw (5) to (6);
		\draw (8) to (9);
		\draw (10) to (11);
		\draw (11) to (12);
		\draw (13) to (14);
		\draw (14) to (15);
  \draw (15) -- (16);
  \draw (18) -- (19);
		\draw (9) to (20);
		\draw (19) to (21);
		\draw (21) to (20);
  \draw (2) to (-2.65,1.5);
  \draw (12) to (-2.65,0);
  \draw (3) to (-1.85,1.5);
  \draw (13) to (-1.85,0);
  \draw (6) to (3.35,1.5);
  \draw (16) to (3.35,0);
  \draw (8) to (4.15,1.5);
  \draw  (18) to (4.15,0);

  \draw (0) [out = 40, in = 140, looseness=0.6] to (20);
  \draw [out = -40, in = -140, looseness=0.6](10) to (21);
\end{tikzpicture}
\caption{$G=C_n \,\square \, K_2$, with vertices labeled as in the proof of Proposition \ref{prop:circularladders}. }
\label{fig:circularladder}
\end{figure}

Without loss of generality, we assume that \poneposs{} first move is not incident to $u_0$ or $v_0$. To describe \ptwo{}'s strategy, we define the following involutive automorphism $\phi: V(G) \to V(G)$. Let $\phi(u_i) = v_{-i}$ for $-k \leq i \leq k$. Note that the only edge fixed by $\phi$ is $u_0v_0$. 

If \pone{}'s previous move was $x \to y$, then \ptwo{} on their turn will move $\phi(x) \to \phi(y)$. Thus \ptwo{} can respond to any move by \pone{} except for if \pone{} moves in either direction along edge $u_0v_0$. Furthermore, \ptwoposs{} first move is not incident to $u_0$ or $v_0$. Suppose without loss of generality \pone{} moves $u_0 \to v_0$ on turn $r$. Since this wasn't \poneposs{} first move, we can assume without loss of generality that \pone{} moved $u_{-1} \to u_0$ on turn $r-2$. On turn $r-1$, \ptwo{} thus moved from $\phi(u_{-1}) = v_{1}$ to $\phi(u_0) = v_0$. Now on turn $r+1$, \ptwo{} moves $v_0 \to v_{-1}$, using up the last available edge incident to $v_0$, as shown in Figure~\ref{fig:copycatladder}. This edge must be unused at this point in the game, since vertex $v_0$ has degree $3$ and was not used as the first vertex in either trail. Then, \pone{} has no possible move on turn $r+2$, and the described strategy is a winning strategy for \ptwo{}. 

 \begin{figure}[ht]
        \centering
            \begin{tikzpicture} [very thick,every node/.style={circle, draw=black, fill=black, inner sep=2}]
		\node [label={[label distance=0.3mm]$u_{-1}$}] (0) at (-4, 1.5) {};
		\node[label={[label distance=1.2mm]$u_{0}$}] (1) at (-2, 1.5) {};
		\node[label={[label distance=1.2mm]$u_{1}$}] (2) at (0, 1.5) {};
		\node[label={[label distance=-1mm]below:$v_{-1}$}] (4) at (-4, 0) {};
		\node[label={[label distance=-0.8mm]below:$v_{0}$}] (5) at (-2, 0) {};
		\node[label={[label distance=-1mm]below:$v_{1}$}] (6) at (0, 0) {};

 		\draw (0) to (4);
  
  \draw[-{Latex[length=3mm]},color=red] (1) -- node[draw=none,fill=none,right]{\pone{}}(5);
		\draw (2) to (6);
  \draw[-{Latex[length=3mm]},color=cyan!70!black] (6) -- node[draw=none,fill=none,above]{\ptwo{}}(5);
		\draw (4) to (5);
  \draw[-{Latex[length=3mm]},color=red] (0) -- node[draw=none,fill=none,above]{\pone{}}(1);
		\draw (2) to (0.5,1.5);
  \draw[-{Latex[length=3mm]},color=cyan!70!black] (5) -- node[draw=none,fill=none,above]{\ptwo{}}(4);
		\draw (1) to (2);
		\draw (6) to (0.5,0);
            \draw (0) to (-4.5,1.5);
            \draw (4) to (-4.5,0);
 
\end{tikzpicture}

\caption{Moves $r-2$ through $r+1$ in the proof of Proposition~\ref{prop:circularladders}, ending with \ptwo{} playing on the last edge incident to $v_0$.}
\label{fig:copycatladder}
\end{figure}

\end{proof}

\section{The Game on Trees}\label{sec tree algorithm}

We now discuss {\sc Trail Trap} on an arbitrary tree $T$. Fraenkel and Simonson \cite{FS93} showed that  {\sc Rooted Partizan Edge Geography} ({\sc Trail Trap} with fixed starting vertices for both players) can be solved in $O(n)$ time on an undirected tree, which trivially provides an $O(n^3)$ time solution to {\sc Trail Trap} on $T$ by iterating over all pairs of starting vertices.
We are able to shave off a factor of $n$ using Lemma~\ref{lem P1 must move towards a center} which shows that if \pone{} has a winning strategy on $T$, then they must first move their token to a central vertex.
Since a tree has at most two central vertices, this provides an $O(n^2)$ time algorithm to resolve {\sc Trail Trap} on $T$ (see Theorem~\ref{thm algorithm for trees}).
To better understand the structure of \pone{}-win trees, we provide a set of necessary conditions for $T$ to be \pone{}-win in Theorem~\ref{thm tree characterization}. We also discuss certain obstacles to a full characterization of \pone{}-win trees.

\subsection{Basic Observations}

A key observation in solving {\sc Trail Trap} on trees is that every trail in a tree is a path, and there is a unique path between each pair of vertices of a tree. The {\em eccentricity} $\ecc_G(v)$ (or simply $\ecc(v)$) of a vertex $v$ in a graph $G$ is $\max\{d(v,x) \colon x \in V(G)\}$. The \textit{radius} of $G$ is $\rad(G) = \min_{v\in V(G)} \ecc(v)$, and the {\em diameter} of $G$ is $\diam(G) = \max_{v \in V(G)} \ecc(v)$.

Note that if $G$ is a tree, then $\MaxTrail(G) = \diam(G)$ and $\MaxTrail(G,v) = \ecc_G(v)$ for any vertex $v$, a fact which we will use extensively later.

The \emph{center} $C(G)$ of a graph $G$ is the set of vertices of minimum eccentricity, i.e., \ $C(G) = \{x \in V(G) \colon \ecc(x) = \rad(x) \}$. A vertex $v \in C(G)$ is called a \emph{central vertex} of $G$. 

\begin{lemma}[Folklore]\label{lem trees have one or two centers}
    For any tree $T$, either (i) $T$ has a single central vertex and $\diam(T) = 2\rad(T)$, or (ii) $T$ has exactly two central vertices, they are adjacent, and $\diam(T) = 2\rad(T) - 1$.
\end{lemma}

In general, {\sc Trail Trap} can be more easily analyzed when \poneposs{} first move occurs on a cut edge. We record a sufficient condition for \ptwo{} to win in the following lemma. 

\begin{lemma}\label{lem P2 and P1 never meet on a graph}
Let $u \to v$ be \poneposs{} first move on a graph $G$. If there is an edge $xy$ in $G - \{uv\}$ such that $v$ and $y$ are in different components of the subgraph $G' = G - \{uv,xy\}$ and $\MaxTrail(G',y) \geq \MaxTrail(G',v)$, then \ptwo{} wins on $G$.
\end{lemma}

\begin{proof}
    Suppose that an edge $xy$ exists in $G - \{uv\}$ such that $v$ and $y$ are in different components of the subgraph $G' = G - \{uv,xy\}$ and such that $\MaxTrail(G',y) \geq \MaxTrail(G',v)$.
    Let $x \to y$ be \ptwoposs{} first move.
    Notice that $ \MaxTrail{(G',v)}$ and $\MaxTrail{(G',y)}$ are, respectively, the maximum number of remaining moves that \pone{} and \ptwo{} can make. Since $v$ and $y$ are in different components of $G'$, neither player can claim an edge from a longest trail available to the other. Thus, \pone{} can make at most $1 + \MaxTrail{(G',v)}$ moves, and \ptwo{} can make at least this many. We conclude that if \poneposs{} first move is $u \to v$, then \ptwo{} wins the resultant partial game on $G$. 
\end{proof}

This lemma is particularly useful for solving {\sc Trail Trap} on graphs with few edges. Recall that, if a graph has at most one edge, then the winner of {\sc Trail Trap} is easily determined, and thus we assume that all graphs have at least two edges.

The following lemma is one of the key conceptual insights into the behavior of {\sc Trail Trap} on a tree: in any winning strategy, \pone{} must start by moving their token to a central vertex, and if there are two central vertices, then \pone{} must move between them.
Figure~\ref{fig:P1 must play into center} depicts two scenarios in which \pone{} does not first move into a central vertex, and one in which they do not first move between two central vertices, leaving \ptwo{} with a winning strategy.

\begin{lemma}\label{lem P1 must move towards a center}
    If $u \to v$ is \poneposs{} first move in {\sc Trail Trap} on a tree $T$ and $v$ is not a central vertex of $T$, then \ptwo{} wins. Furthermore, if $T$ has two central vertices $c$ and $c'$ and \pone{} plays neither $c\to c'$ nor $c'\to c$ as their first move, then \ptwo{} wins. 
\end{lemma}
\begin{proof}
    Our goal is to apply Lemma \ref{lem P2 and P1 never meet on a graph}, so we must show that if $u\to v$ is an edge in $T$ such that $v$ is not a central vertex, there is an edge $ab\in T$ such that $v$ and $b$ are in different components of $T - \{uv, ab\}$ and $\MaxTrail(T - \{uv, ab\}, v) \leq \MaxTrail(T - \{uv, ab\}, b)$. 
    
    Let $c$ be the central vertex of $T$ nearest to $v$. Let $r$ be the radius of $T$, and let $k = d_T(c,v)$. By assumption, $k > 0$ as $v$ is not a central vertex. Let $vPc$ be the unique path from $v$ to $c$. 
    
    If $u\notin vPc$, let $vx$ be the edge incident to $v$ in $vPc$. See Figure \ref{subfigA:P1 center} for an illustration of this case. We claim that $vx$ is our desired edge for Lemma \ref{lem P2 and P1 never meet on a graph}. First observe that $x$ is not in the same component of $T - \{uv,  vx\}$ as $v$. Second, if $T'$ is the component of $T - \{uv, vx\}$ which contains $v$, we observe that $\MaxTrail(T',v) = \ecc_{T'}(v) \leq r - k$ as every path $wPc$ in $T$ from a vertex $w\in T'$ to $c$ must contain $v$, and every such path has length at most $r$. By definition of a central vertex, there are at least two edge-disjoint paths of length $r$ which end at $c$ if $c$ is the unique central vertex, or two edge-disjoint paths of length $r - 1$ which end at each of the central vertices and do not use the edge between the two central vertices. In either case, there is a path $P'$ of length $r$ which ends at $c$ and is edge-disjoint from $vPc$. If $T''$ is the component of $T - \{uv, vx\}$ containing $x$, it follows that $\MaxTrail(T'',x) \geq |xPc| + |P'| = k - 1 + r$. As $k > 0$, it follows that $\MaxTrail(T'',x) \geq r - 1 + k > r - k \geq \MaxTrail(T',v)$. By Lemma \ref{lem P2 and P1 never meet on a graph}, \ptwo{} wins. 
    
    Now assume that $u\in vPc$. See Figure \ref{subfigB:P1 center} for an illustration of this case. Let $ux$ be the first edge in the path $uPc$. The vertex $x$ is not in the same component of $T - \{uv, ux\}$ as $v$, and by identical reasoning as before, $\MaxTrail(T - \{uv,ux\},x) \geq k - 2 + r  \geq r - k \geq  \MaxTrail(T - \{uv,ux\},v)$. Hence, we can apply Lemma \ref{lem P2 and P1 never meet on a graph} to conclude that \ptwo{} also wins in this scenario.

    For the second claim, assume that $T$ has two central vertices $c$ and $c'$, and that \poneposs{} first move is an edge $u\to c$ where $u \neq c'$. See Figure \ref{subfigC:P1 center} for an illustration of this case. We claim that \ptwo{} can win by first moving $c \to c'$. By definition of a central vertex, there are at least two edge-disjoint paths
    of length $r - 1$ in $T$, which end at $c$ and $c'$ and do not use the edge $cc'$. As all paths from $c$ which do not contain $cc'$ have length at most $r - 1$, it follows that $\MaxTrail(T - \{uc, cc'\}, c) \leq r - 1$. Since $\MaxTrail(T - \{uc, cc'\},c') = r - 1$, \ptwo{} wins by Lemma \ref{lem P2 and P1 never meet on a graph}.
\end{proof}

\begin{figure}
    \centering
    \begin{subfigure}{.3\textwidth}
        \centering
        \begin{tikzpicture}
        [very thick,every node/.style={circle, draw=black, fill=black, inner sep=0, minimum size=3pt}, scale=.8]

        \node [label=below:{$c$ \ }] (c) at (0,0) {};
        \node (d0) at (5,3) {};
        \node (d1) at (5,1.32) {};
        \node (d2) at (5,.6) {};
        \node (d3) at (5,0) {};
        \node (d4) at (5,-2) {};
        \node [label=below:{$u$}] (u) at (4,0) {};
        \node [label=below:{$v$}] (v) at (2.8,0) {};

        \node[label=below:{$x$}] (x) at (1.6,0) {};

        \draw (d0) -- (c) -- (d4);
        \draw (c) -- (1.6, 0);
        \draw (u) -- (d3);
        \draw (v) -- (d1);
        \draw (u) -- (d2);

        \draw[ultra thick][-{Latex[length=3mm]},color=red] (u) -- node[draw=none,fill=none,below]{\pone{}}(v);
        \draw[ultra thick][-{Latex[length=3mm]},color=cyan!70!black] (v) -- node[draw=none,fill=none,above]{\ptwo{}} (x);

        \node[fill=none, draw=none] (r) at (2.55, -1.45) {$r$};
        \end{tikzpicture}
        \caption{$u \notin vPc$}
        \label{subfigA:P1 center}
    \end{subfigure}
    \hfill
    \begin{subfigure}{.3\textwidth}
        \centering
        \begin{tikzpicture}
        [very thick,every node/.style={circle, draw=black, fill=black, inner sep=0, minimum size=3pt}, scale=.8]

        \node [label=below:{$c$ \ }] (c) at (0,0) {};
        \node (d0) at (5,3) {};
        \node (d1) at (5,1.32) {};
        \node (d2) at (5,.6) {};
        \node (d3) at (5,0) {};
        \node (d4) at (5,-2) {};
        \node [label=below:{$u$}] (u) at (2.8,0) {};
        \node [label=below:{$v$}] (v) at (4,0) {};

        \node[label=below:{$x$}] (x) at (1.6,0) {};

        \draw (d0) -- (c) -- (d4);
        \draw (c) -- (1.6,0);
        \draw (v) -- (d3);
        \draw (u) -- (d1);
        \draw (v) -- (d2);

        \draw[ultra thick][-{Latex[length=3mm]},color=red] (u) -- node[draw=none,fill=none,below]{\pone{}}(v);
        \draw[ultra thick][-{Latex[length=3mm]},color=cyan!70!black] (u) -- node[draw=none,fill=none,above]{\ptwo{}} (x);

        \node[fill=none, draw=none] (r) at (2.55, -1.45) {$r$};
        
        \end{tikzpicture}
        \caption{$u \in vPc$}
        \label{subfigB:P1 center}
    \end{subfigure}
    \hfill
    \begin{subfigure}{.3\textwidth}
        \centering
        \begin{tikzpicture}
        [very thick,every node/.style={circle, draw=black, fill=black, inner sep=0, minimum size=3pt}, scale=.8]

        \node [label=above:{$c$ \ }] (c) at (0,1.2) {};
        \node [label=below:{$c'$ \ }] (c1) at (0,-.2) {};
        \node (d0) at (5,3) {};
        \node (d1) at (5,1.2) {};
        \node (d2) at (5,-2) {};
        \node [label=above:{$u$}] (u) at (1.2,1.632) {};

        \node[fill=none, draw=none] (r) at (2.5, .9) {$r-1$};
        \node[fill=none, draw=none] (r) at (2.5, -1.6) {$r-1$};

        \draw (d2) -- (c1);
        \draw (c) -- (d1);
        \draw (u) -- (d0);

        \draw[ultra thick][-{Latex[length=3mm]},color=red] (u) -- node[draw=none,fill=none,above]{\pone{}}(c);
        \draw[ultra thick][-{Latex[length=3mm]},color=cyan!70!black] (c) -- node[draw=none,fill=none,left]{\ptwo{}} (c1);
        
        \end{tikzpicture}
        \caption{$v = c$, center is $\{c, c'\}$}
        \label{subfigC:P1 center}
    \end{subfigure}
    \caption{The scenarios considered in the proof of Lemma~\ref{lem P1 must move towards a center}. If \pone{} does not play into a central vertex $c$, or does not play between two central vertices $c,c'$, then \ptwo{} wins.}
    \label{fig:P1 must play into center}
\end{figure}

With significant care, one can describe broader classes of graphs $G$ for which \ptwo{} wins whenever \pone{} does not start by playing into a central vertex.
While it would be interesting to determine which \pone{}-win graphs require \pone{} to play into a central vertex on their first move, a complete characterization of such graphs appears a complicated task.

We now recall Proposition~\ref{lem:p1 loses on degree 1 and 2}: if $u \to v$ is \poneposs{} first move in a winning strategy on a graph $G$, then $\deg(v) \geq 3$.
By Lemma~\ref{lem P1 must move towards a center}, if $G$ is a \pone{}-win tree, then $v$ must also be a central vertex.
We now prove additional conditions on the degrees of the central vertices in a \pone{}-win tree.

\begin{lemma}\label{lemma central vertices general}
    Let $u\to v$ be \poneposs{} first move on a tree $T$ with at least three vertices. If $v$ is not a central vertex of degree $3$, then \ptwo{} wins. Furthermore, if $u$ and $v$ are both central vertices of $T$ and $\deg(u) \neq 2$, then \ptwo{} again wins.
\end{lemma}
\begin{proof}
    By Lemma~\ref{lem P1 must move towards a center}, we may assume that $v$ is a central vertex. If $\deg(v) \leq 2$, then \ptwo{} wins by Proposition~\ref{lem:p1 loses on degree 1 and 2}. If $\deg(v) > 3$, then we have two cases to consider. If there exists two paths $P,P'$ in $T - \{uv\}$ which start at $v$ and have length $\MaxTrail(T - \{uv\},v)$, let $w$ be a neighbor of $v$ not on either path. By construction, \ptwo{} may play a trail of length $\MaxTrail(T - \{uv\}, v) + 1$ starting with the edge $wv$  whereas \pone{} may at best play a trail of total length $\MaxTrail(T - \{uv\}, v) + 1$. Thus \ptwo{} wins. If there is only one path $P$ starting at $v$ of length $\MaxTrail(T - \{uv\},v)$, let $w$ be the neighbor of $v$ on the path $P$. Then $uv,vw$ is a pair of edges which satisfies the conditions of Lemma \ref{lem P2 and P1 never meet on a graph}. Indeed, $v$ and $w$ are in different connected components of $T - \{uv, vw\}$, and $\MaxTrail(T - \{uv, vw\},w)  \geq \MaxTrail(T - \{uv, vw\},v)$. Thus \ptwo{} again wins. 

    For the second claim, recall that if $T$ has radius $r$ and two central vertices $u$ and $v$, then $\MaxTrail(T - \{uv\},u) = \MaxTrail(T - \{uv\},v) = r - 1$. Let $P_u$ be a path of length $r - 1$ in $T - \{uv\}$ starting at $u$, and let $w$ be a neighbor of $u$ not in $P_u \cup \{v\}$. Such a vertex exists since $T \neq K_2$ and $\deg(u) \neq 2$ implies $\deg(u) > 2$. We claim that the edges $uv$ and $wu$ satisfy the conditions of Lemma \ref{lem P2 and P1 never meet on a graph}. Indeed, $u$ and $v$ are in different components of $T - \{uv, wu\}$, and $\MaxTrail(T - \{uv, wu\},u) = r - 1 = \MaxTrail(T - \{uv, wu\},v) $. Thus by Lemma \ref{lem P2 and P1 never meet on a graph}, \ptwo{} wins on $T$.
\end{proof}

\begin{figure}
    \centering
    \begin{tikzpicture}
        [very thick,every node/.style={circle, draw=black, fill=black, inner sep=0, minimum size=3pt}, scale=.8]

        \node (c) at (0,0) {};
        \node (0) at (2,2) {};
        \node (1) at (2,0) {};
        \node (2) at (2,-2) {};
        \node (00) at (4,2.5) {};
        \node (01) at (4,1.5) {};
        \node (10) at (4,.5) {};
        \node (11) at (4,-.5) {};
        \node (20) at (4,-1.5) {};
        \node (21) at (4,-2.5) {};
        \node (000) at (6,2.5) {};
        \node (010) at (6,1.5) {};
        \node (100) at (6,.5) {};
        \node (110) at (6,-.5) {};
        \node (200) at (6,-1.5) {};
        \node (210) at (6,-2.5) {};

        \foreach \i in {0,1,2}
        {
        \draw (c) -- (\i);
        \draw (\i) -- (\i 0) -- (\i 00);
        \draw (\i) -- (\i 1) -- (\i 10);
        }

        \draw[ultra thick][-{Latex[length=3mm]},color=red] (0) -- node[draw=none,fill=none,above]{\pone{} \ }(c);

        \draw[ultra thick][-{Latex[length=3mm]},color=cyan!70!black] (000) -- node[draw=none,fill=none,above]{\ptwo{}} (00);

        \end{tikzpicture}
    \caption{A tree with a single central vertex which does not meet condition {\em (i)} to be \pone{}-win in Theorem~\ref{thm tree characterization}.}
    \label{fig:behind_the_back}
\end{figure}

We now give a number of necessary conditions for a tree to be \pone{}-win.

\begin{theorem}\label{thm tree characterization}
    Let $T$ be a \pone{}-win tree with radius $r$.
    If $T$ has a single central vertex $c$, then $\deg(c) = 3$. Further, for some neighbor $x$ of $c$, 
    \begin{enumerate}[(i)]
        \item the component of $T - \{cx\}$ containing $x$ has diameter at most $r$, and
        \item 
        the component of $T - \{cx\}$ containing $c$ has just one central vertex, namely $c$.
    \end{enumerate}
    If $T$ has two central vertices $c_1$ and $c_2$ with $\deg(c_1) \geq \deg(c_2)$, then $\deg(c_1) = 3$ and $\deg(c_2) = 2$. Further, 
    \begin{enumerate}[(i)]
        \item[(iii)] the component of $T - \{c_1c_2\}$ containing $c_2$ has diameter at most $r-1$, and
        \item[(iv)] the component of $T-\{c_1c_2\}$ containing $c_1$ has just one central vertex, namely $c_1$. 
    \end{enumerate}
    
\end{theorem}

\begin{proof}
    By Lemmas~\ref{lem P1 must move towards a center} and~\ref{lemma central vertices general}, in order for a tree $T$ to be \pone{}-win, either $T$ has a single central vertex $c$ of degree $3$ and \poneposs{} first move ends at $c$, or $T$ has two central vertices $c_1$ and $c_2$ of degrees $3$ and $2$, respectively, and \poneposs{} first move is $c_2 \to c_1$.

    Suppose that $T$ has a single central vertex $c$ with $N(c) = \{x_1, x_2, x_3\}$.
    If, for all $i$, the component of $T - \{cx_i\}$ that contains $x_i$ has diameter strictly greater than $r$, then no matter which move $x_i \to c$  \pone{} makes first, \ptwo{} can construct a trail of length $r + 1$ in the component of $T-\{cx_i\}$ containing $x_i$. As \pone{} must then play in the component containing $c$, \pone{} can construct a trail of length at most $r + 1$ in total (including the edge $x_i\to c$), so \pone{} loses in this case. 
    Without loss of generality, suppose the component of $T - \{cx_1\}$ containing $x_1$ has diameter at most $r$.
    If $c$ is not the unique central vertex of its component of $T - \{cx_1\}$, then exactly one of $x_2$ or $x_3$ has eccentricity $r-1$ in $T-c$, say $x_2$. \ptwoposs{} winning strategy is to start with $c \to x_2$. This means that \pone{} can construct a trail of length at most $r$, while \ptwo{} can construct a trail of length at least $r$.

    Suppose that $T$ has two central vertices $c_1, c_2$ with $\deg(c_1) = 3$ and $\deg(c_2) = 2$.
    Recall that $c_2 \to c_1$ is \poneposs{} first move in any winning strategy.
    If the diameter of the component of $T - \{c_1c_2\}$ containing $c_2$ is at least $r$, then \ptwo{} can find a trail of length at least $r$ in this component, and \pone{} has at most $r - 1$ moves remaining in the component containing $c_1$.
    Furthermore, if $c_1$ is not the unique central vertex of its component in $T - \{c_1c_2\}$, then \ptwo{} has a winning strategy, as described above.
\end{proof}

Theorem~\ref{thm tree characterization} provides a set of necessary conditions for a tree $T$ to be \pone{}-win.
Figure~\ref{fig:behind_the_back} depicts an example of a tree with a single central vertex which does not meet condition {\em (i)} in Theorem~\ref{thm tree characterization}, and thus is \ptwo{}-win.

The set of \pone{}-win trees does not seem to have a concise characterization. 
For instance, if a tree $T$ with a single central vertex $c$ meets the necessary conditions of Theorem~\ref{thm tree characterization}, and if $c$ is contained in every path of length greater than $r$ in its component of $T - \{cx\}$, then $T$ is \pone{}-win.
On the other hand, if the latter condition does not hold, then
to determine the winner on $T$ we must consider the case in which \ptwo{} plays elsewhere in the tree (specifically, on the path of length $>r$ which does not contain $c$).
If \pone{} does not move their token towards \ptwoposs{} on their second move, then \ptwo{} has a winning strategy for the resulting partial game. 
If \pone{} does move their token towards \ptwo{}, then further analysis is required to analyze the interactions between the players. 

\begin{figure}[hbt]
    \centering
    \begin{subfigure}{.45\textwidth}
        \centering
        \begin{tikzpicture}
        [very thick,every node/.style={circle, draw=black, fill=black, inner sep=0, minimum size=3pt}, scale=.8]

        \node[label={$c$}] (c) at (0,0) {};
        \node (0) at (2,1.5) {};
        \node (1) at (2,0) {};
        \node (2) at (2,-1.5) {};
        \node (00) at (3.33,2) {};
        \node (01) at (3.33,1) {};
        \node (20) at (3.33,-1.5) {};
        \node (000) at (4.66,2) {};
        \node (010) at (4.66,1) {};
        \node (200) at (4.66,-1.5) {};

        \draw (000) -- (00) -- (0) -- (c) -- (2) -- (20) -- (200);
        \draw (0) -- (01) -- (010);
        \draw (1) -- (c);

        \draw[ultra thick][-{Latex[length=3mm]},color=red] (1) -- node[draw=none,fill=none,above]{\pone{}}(c);
        
        \end{tikzpicture}
        
        \caption{A \pone{}-win tree}
        \label{subfig:tree characterization p1win}
    \end{subfigure}
    \hfill
    \begin{subfigure}{.45\textwidth}
        \centering 
        \begin{tikzpicture}
        [very thick,every node/.style={circle, draw=black, fill=black, inner sep=0, minimum size=3pt}, scale=.8]

        \node[label={$c$}] (c) at (0,0) {};
        \node (0) at (2,1.5) {};
        \node (1) at (2,0) {};
        \node (2) at (2,-1.5) {};
        \node (00) at (3.33,2) {};
        \node (01) at (3.33,1) {};
        \node (20) at (3.33,-1.5) {};
        \node (200) at (4.66,-1.5) {};
        \node (2000) at (6, -1.5) {};
        \node (20000) at (7.33,-1.5) {};

        \node (000) at (4.66,2) {};
        \node (0000) at (6,2) {};
        \node (00000) at (7.33,2) {};
        \node (010) at (4.66,1) {};
      
        \foreach \i in {0,2}
        {
        \draw (c) -- (\i);
        \draw (\i) -- (\i 0) -- (\i 00) -- (\i 000) -- (\i 0000);
        }

        \draw (0) -- (01) -- (010);

        \draw (1) -- (c);

        \draw[ultra thick][-{Latex[length=3mm]},color=red] (1) -- node[draw=none,fill=none,above]{\pone{}}(c);

        \draw[ultra thick][-{Latex[length=3mm]},color=cyan!70!black] (00000) -- node[draw=none,fill=none,above]{\ptwo{}} (0000);
        
        \end{tikzpicture}
        \caption{A \ptwo{}-win tree}
        \label{subfig:tree characterization p2win}
    \end{subfigure}
    \caption{Two trees meeting the necessary (but not sufficient) conditions of Theorem~\ref{thm tree characterization} to be \pone{}-win.}
    \label{fig:sufficiency_tree_thm}
\end{figure}

Figure~\ref{fig:sufficiency_tree_thm} depicts two such scenarios.
Let $T$ be the tree in Figure~\ref{subfig:tree characterization p1win}. To determine that $T$ is \pone{}-win, we must consider the case in which the two players move towards one another on the top branch from $c$: if \ptwo{} starts by moving from a leaf along the path of length $4$ in this branch, then \pone{} must move towards \ptwo{} to win. 
Having done so, \pone{} is able to stop \ptwo{} from taking this path, and wins.
In Figure~\ref{subfig:tree characterization p2win}, while the tree does meet the necessary conditions of Theorem~\ref{thm tree characterization}, parity conditions on the path from \ptwoposs{} first move to the central vertex determine the winner.
We leave this as an exercise for the reader.

We now turn our attention to the algorithmic problem of determining the winner of {\sc Trail Trap} on a tree. It suffices, by 
Lemma~\ref{lem P1 must move towards a center}, to consider the case where \pone{} starts by moving their token to 
a central vertex. 
Once \ptwo{} makes their first move, determining the winner of {\sc Trail Trap} is equivalent to determining the winner of a game of {\sc Rooted Partizan Edge Geography}. As part of providing an algorithm for a Geography variant on directed graphs, Fraenkel and Simonson \cite{FS93} give an $O(n)$ time algorithm for determining the winner of {\sc Rooted Partizan Edge Geography} on an undirected tree. Hence, we have an efficient algorithm to resolve the trees left unresolved by the partial characterization of \pone{}-win trees in Theorem \ref{thm tree characterization}.
One can naively determine the winner of {\sc Trail Trap} on a tree in $O(n^3)$ time by running Fraenkel and Simonson's algorithm once for each pair of starting moves. 
However, the center of a tree can be determined in linear time~{\cite{CHH1981}}, allowing us to only vary over the starting moves for \ptwo{} and shave off a factor of $n$.
        \begin{theorem}\label{thm algorithm for trees}
            The winner of {\sc Trail Trap} on a tree $T$ can be determined in $O(n^2)$ time. 
        \end{theorem}
        \begin{proof}
            
            We informally sketch the algorithm.  We assume that we have a $O(n)$-time subroutine $RPEG$ which takes in a tree $T$, \poneposs{} first move $u \to v$, and \ptwoposs{} first move $x \to y$ and returns \pone{}-win or \ptwo{}-win. Note that $u,v,x,y$ may not be distinct. We have that $RPEG(T,u\to v,x\to y)$ returns \pone{}-win if and only if \pone{} wins  {\sc Rooted Partizan Edge Geography} on $T$ with the given starting moves $u\to v$ and $x\to y$. 
            
            Our algorithm first computes the central vertices of $T$.  This can be done in $O(n)$ time \cite{CHH1981}.
            
             If $T$ has one central vertex $c$, our algorithm first checks if $c$ has degree $3$, and returns \ptwo{}-win if not. If $c$ has degree $3$, our algorithm iterates over \poneposs{} starting moves $u \to c$ for each neighbor $u$ of $c$, then iterates over all starting moves $x\to y$ for \ptwo{}, and finally runs $RPEG(T,u\to c, x\to y)$. If there is a starting move $x\to y$ for \ptwo{} such that $RPEG(T,u\to c, x\to y)$ returns \ptwo{}-win, then our algorithm considers a different starting edge for \pone{}. If \ptwo{} has a winning starting edge for each of the three choices of a starting edge for \pone{}, then our algorithm returns \ptwo{}-win. Otherwise, our algorithm returns \pone{}-win. 
            
            As there are three possible starting edges for \pone{} and $O(n)$ starting edges for \ptwo{} (as $T$ is a tree), our algorithm runs in $O(n^2)$ time on a tree with one central vertex. 

            If $T$ has two central vertices $c_1$ and $c_2$ with $\deg(c_1)\geq \deg(c_2)$, our algorithm checks if $\deg(c_1) = 3$ and $\deg(c_2) = 2$, and returns \ptwo{}-win if not. If these degree conditions are satisfied, the algorithm then iterates over all starting moves $x\to y$ for \ptwo{} and runs $RPEG(T,c_2\to c_1, x\to y)$. If there is an edge $x\to y$ such that $RPEG(T, c_2 \to c_1, x \to y)$  returns \ptwo{}-win, then our algorithm returns \ptwo{}-win and returns \pone{}-win otherwise.
            
            As there is exactly one possible starting move for \pone{} by Lemma \ref{lem P1 must move towards a center} and $O(n)$ possible starting moves for \ptwo{}, our algorithm runs in $O(n^2)$ time on a tree with two central vertices. 

            The correctness of the algorithm then follows from  Lemmas \ref{lem P1 must move towards a center} and \ref{lemma central vertices general}. 
        \end{proof}

\section{NP-hardness} \label{sec NP-hardness}

In this section, we prove that deciding whether a graph is \ptwo{}-win is NP-hard, even for connected bipartite planar graphs with maximum degree $4$.
We will consider the following decision problems. Recall that a Hamiltonian path in a graph is a path that visits every vertex of the graph exactly once.

\medskip
	\begin{tabular}{l}
		\textsc{Hamiltonian path}\\
		Input: A graph $G$.\\ \vspace{0.5cm}
		Question: Does $G$ have a Hamiltonian path?\\
        

        \textsc{Hamiltonian cycle through specified edge}\\
            Input: A graph $G$ and edge $e$ in $G$.\\ \vspace{0.5cm}
            Question: Does $G$ contain a Hamiltonian cycle through $e$?\\

        \textsc{Fixed long trail}\\
		Input: A graph $G$ and a vertex $u \in V(G)$.\\ \vspace{0.5cm}
		Question: Does $G$ have a trail of length at least $|V(G)| + 1$ starting at $u$?\\

        \textsc{Trail Trap \ptwo{}-win}\\
		Input: A graph $G$.\\
		Question: Does \ptwo{} win {\sc Trail Trap} on $G$?\\

        \medskip

	\end{tabular}


Before we prove that {\sc Trail Trap \ptwo{}-win} is NP-hard, we note two results from the literature which establish the computational complexity of {\sc Hamiltonian path} and {\sc Fixed long trail} for planar cubic bipartite graphs. 
We recall that if a problem is NP-complete, then it is also NP-hard.

In a comment on the Theoretical Computer Science Stack Exchange, Labarre proved the following result.

\begin{theorem}[\cite{Labarre}]\label{thm:labarre}
    \textsc{Hamiltonian cycle through specified edge} is NP-complete for planar cubic bipartite graphs.
\end{theorem}


Munaro~\cite[Theorem~23]{MUNARO20171210} used Theorem \ref{thm:labarre} to prove that {\sc Hamiltonian path} is NP-complete for planar cubic bipartite graphs, by replacing an edge $uv$ in a planar cubic bipartite graph $G$ with the gadget depicted in Figure~\ref{fig:fixedlongtrail}.
Below, we state the key idea of Munaro's proof as a lemma. 

\begin{lemma}[{\cite[p.~1221]{MUNARO20171210}}]\label{lem:munaro}
    Let $uv$ be an edge in a planar cubic bipartite graph $G$, and let $G'$ be the graph obtained by replacing $uv$ with the gadget depicted in Figure~\ref{fig:fixedlongtrail}.
    The graph $G$ contains a Hamiltonian cycle through $uv$ if and only if $G'$ contains a Hamiltonian path.
\end{lemma}

We are able to use this lemma to show that {\sc Fixed long trail} is NP-complete for the same class of graphs.

\begin{figure}[hbt]
    \centering
    \begin{tikzpicture}
        [thick,every node/.style={circle, draw=black, fill=black, inner sep=0, minimum size=3pt}, hl/.style={gray!50, line width=5pt}, scale=.75]


        \coordinate (A) at (-6.5,0);
        \coordinate (B) at (0,1.5);
        \coordinate (C) at (6.5,0);
        \coordinate (D) at (0,-1.5);

        \foreach \p in {A, B, C, D}
        {
        \node (g1\p) at ($(\p) + (-2,0)$) {};
        \node (g2\p) at ($(\p) + (-1,1)$) {};
        \node (g3\p) at ($(\p) + (1,1)$) {};
        \node (g4\p) at ($(\p) + (2,0)$) {};
        \node (g5\p) at ($(\p) + (1,-1)$) {};
        \node (g6\p) at ($(\p) + (-1,-1)$) {};
        \node (g7\p) at ($(\p) + (-.5,0)$) {};
        \node (g8\p) at ($(\p) + (.5,0)$) {};
        }

        \node[fill=none, minimum size=4pt] (x) at ($(A) + (3.5,0)$) {};
        \node[draw=none,fill=none] (xlabel) at ($(x) + (-.15,.25)$) {$x$};
        \node (y) at ($(C) + (-3.5,0)$) {};
        \node[draw=none, fill=none] (ylabel) at ($(y) + (.15,.35)$) {$y$};

        \node[draw=none,fill=none] (start) at ($(x) + (.125,.1875)$) {};
        \node[draw=none,fill=none] (end) at ($(y) + (.2,0)$) {};

        \node[draw=none,fill=none] (x'label) at ($(g1B) + (-.15,.25)$) {$x'$};
        \node[draw=none,fill=none] (y'label) at ($(g1C) + (0,.4)$) {$y'$};

        \node[label={$u$}] (u) at ($(A) + (-3,0)$) {};
        \node[label={$v$}] (v) at ($(C) + (3,0)$) {};

        \draw[hl][{Bar[width=4mm,length=0mm]}-{Latex[length=4mm]}] (start) -- (g1B);
    
        \draw[hl][-{Latex[length=4mm] . Bracket[width=4mm, length=0mm]}] 
        (g1C) -- (end);
        
        \draw[hl] 
        (g1B) -- (g2B) -- (g7B) -- (g6B) -- (g5B) -- (g8B) -- (g3B) -- (g4B) -- (y) -- (g4D) -- (g5D) -- (g8D) -- (g3D) -- (g2D) -- (g7D) -- (g6D) -- (g1D) -- (x) -- (g4A) -- (g3A) -- (g8A) -- (g5A) -- (g6A) -- (g7A) -- (g2A) -- (g1A) -- (u);


        
        \draw[hl][-{Latex[length=4mm]}] 
        (u) -- ($(u) + (-1,0)$);
        \draw[hl][-{Latex[length=4mm]}]  
        ($(v) + (1,0)$) -- (v);
        
        \draw[hl] 
        (v) -- (g4C) -- (g3C) -- (g8C) -- (g5C) -- (g6C) -- (g7C) -- (g2C) -- (g1C);

        \foreach \p in {A, B, C, D}
        {
        \draw (g1\p) -- (g2\p) -- (g3\p) -- (g4\p) -- (g5\p) -- (g8\p) -- (g7\p) -- (g6\p) -- (g1\p);
        \draw (g2\p) -- (g7\p);
        \draw (g3\p) -- (g8\p);
        \draw (g5\p) -- (g6\p);
        }

        \draw (g4A) -- (x) -- (g1B);
        \draw (x) -- (g1D);
        \draw (g1C) -- (y) -- (g4B);
        \draw (y) -- (g4D);

        \draw (u) -- (g1A);
        \draw ($(A) + (-3.75,.5)$) -- (u) -- ($(A) + (-3.75,-.5)$);
        \draw (v) -- (g4C);
        \draw ($(C) + (3.75,.5)$) -- (v) -- ($(C) + (3.75,-.5)$);

        
    \end{tikzpicture}
    \caption{The gadget used in~\cite[Theorem 23]{MUNARO20171210} to replace the edge $uv$ in a planar cubic bipartite graph. A fixed long trail from $x$ (which first takes the edge $xx'$, exits the gadget at $u$, follows a Hamiltonian path in $G$ from $u$ to $v$, and ends with the edge $y'y$) is highlighted (see proof of Theorem~\ref{thm:fixedlongtrail}).}
    \label{fig:fixedlongtrail}
\end{figure}


\begin{theorem}\label{thm:fixedlongtrail}
    Let $uv$ be an edge in a planar cubic bipartite graph $G$, and let $G'$ be the planar cubic bipartite graph obtained from $G$ by replacing the edge $uv$ with the gadget depicted in Figure~\ref{fig:fixedlongtrail}.
    The graph $G$ contains a Hamiltonian cycle through $uv$ if and only if $G'$ contains a trail of length $|V(G')| + 1$
    starting from the vertex $x$ (also depicted in Figure~\ref{fig:fixedlongtrail}).
\end{theorem}

\begin{proof}
    First, suppose that $G$ contains a Hamiltonian cycle through $uv$.
    Deleting the edge $uv$ from the cycle produces a Hamiltonian path $P$ in $G$ with endpoints $u$ and $v$.
    The path $P$ can be extended to a long trail of length $|V(G')|+1$ in $G'$ with endpoint $x$ by following the highlighted path in Figure~\ref{fig:fixedlongtrail}.

    For the converse, suppose $G'$ contains a trail $T$ of length $|V(G')|+1$ starting from $x$.
    Let $z$ denote the other endpoint of $T$.
    Note that the only vertices of degree $3$ in $T$ are the endpoints of the trail.
    Let $x'$ be the second vertex visited by $T$, and let $z'$ be the penultimate vertex visited by $T$.
    Then $T - \{xx', zz'\}$ is a Hamiltonian path in $G'$.
    By Lemma~\ref{lem:munaro}, 
    $G$ contains a Hamiltonian cycle through the edge $uv$, which completes the proof.
\end{proof}

\begin{corollary}
    \label{cor:npc-trail-cubic-fixed-start}
    \textsc{Fixed long trail} is NP-complete on bipartite cubic planar graphs.
\end{corollary}

\begin{proof}
  Given a trail in a graph $G$, we can check in polynomial time if the trail starts at $x$ and has length at least $|V(G)|+1$. Thus, \textsc{Fixed long trail} is in NP.

  Let $G$ be a bipartite cubic planar graph with a fixed edge $uv$, and let $G'$ be the bipartite cubic planar graph formed from $G$ by replacing the edge $uv$ with the gadget shown in Figure~\ref{fig:fixedlongtrail}.
  By Theorem~\ref{thm:fixedlongtrail}, $G$ contains a Hamiltonian cycle through the fixed edge $uv$ if and only if $G'$ contains a fixed long trail starting from the vertex $x$.
  Since determining if $G$ has a Hamiltonian cycle through $uv$ is NP-hard by Theorem~\ref{thm:labarre}, determining whether $G'$ has a trail starting at $x$ of length $|V(G')|+1$ is also NP-hard.
  Hence, \textsc{Fixed long trail} is NP-complete.
\end{proof}

We are now able to prove that {\sc Trail Trap \ptwo{}-win} is NP-hard on connected bipartite planar graphs with maximum degree 4. Note that this refers to the more restrictive setting of whether \ptwo{} wins from the starting position (i.e. a graph with no player tokens yet) rather than to partial games, as is common elsewhere in the literature.

\begin{theorem}
    \label{thm:NP}
    {\sc Trail Trap \ptwo{}-win} is NP-hard  on connected bipartite planar graphs with maximum degree 4.
\end{theorem}
\begin{proof}
    Let $G$ be a connected bipartite planar cubic graph on $n$ vertices and fix a vertex $w \in V(G)$. We construct the graph $\widehat{G}$ from $G$ by adding a new vertex $u$ adjacent to $w$, i.e.~$V(\widehat{G}) = V(G)\cup\{u\}$ and $E(\widehat{G}) = E(G)\cup\{uw\}$. Since $G$ is cubic, the maximum trail length of $G$ is at most $n+1$, so every trail in $\widehat{G}$ that does not start or end at $u$ is of length at most $n+1$ (since $\deg_{\widehat{G}}(u) = 1$). Therefore $\widehat{G}$ has a trail of length $n+2$ if and only if $G$ has a trail of length $n+1$ starting at $w$. By Corollary~\ref{cor:npc-trail-cubic-fixed-start}, determining whether $\widehat{G}$ has a trail of length $n+2$ is NP-hard, since it is equivalent to determining whether $G$ has a trail of length $n+1$ starting at $w$.

    We reduce finding a trail of length $n+2$ in $\widehat{G}$ to {\sc Trail Trap \ptwo{}-win} on some connected bipartite planar graph (with maximum degree 4)  as follows. Let $Q$ be the path with vertices $v_1, \ldots, v_{2n+3}$ and edges $v_i v_{i+1}$ for $1 \le i \le 2i+2$. Since $v_{n+2}$ is the unique central vertex of $Q$, we refer to it as $c$. Let $G'$ be the graph obtained from the disjoint union of $Q$ and $\widehat{G}$ by adding an edge $u c$ (see Figure \ref{fig:G'}). Observe that $G'$ is connected, bipartite, planar and has maximum degree 4. We claim that \ptwo{} wins {\sc Trail Trap} on $G'$ if and only if $\widehat{G}$ has a trail of length $n+2$.

\begin{figure}[h!!]
\centering
\begin{tikzpicture}
[
vert/.style={circle,fill=black,draw=black, inner sep=0.05cm}
] 

\filldraw[color=black, fill=black!5!white, thick] (5,-3) circle (2);
\filldraw[color=black, fill=black!15!white, thick] (5,-3.5) circle (1.25);

\node[vert, label={90:$v_1$}] (v1) at (0,0) {};
\node[vert, label={90:$v_2$}] (v2) at (1,0) {};
\node (v3) at (2,0) {$\cdots$};
\node[vert, label={90:$v_{n+1}$}] (v4) at (4,0) {};
\node[vert, label={90:$c$}] (v5) at (5,0) {};
\node[vert, label={90:$v_{n+3}$}] (v6) at (6,0) {};
\node (v7) at (8,0) {$\cdots$};
\node[vert, label={90:$v_{2n+2}$}] (v8) at (9,0) {};
\node[vert, label={90:$v_{2n+3}$}] (v9) at (10,0) {};
\node(Q) at (8.5,-0.45) {$Q$};

\foreach \x [remember=\x as \lastx (initially 1)] in {1,...,9}
\path (v\x) edge (v\lastx);

\node[vert, label={0:$u$}] (u) at (5,-1.5) {};
\node[vert, label={0:$w$}] (w) at (5,-2.5) {};
\draw (u) -- (w);
\node (g) at (4.5,-1.75) {$\widehat{G}$};
\node (g) at (5,-3.75) {$G$};

\path (u) edge (v5);

\end{tikzpicture}

\caption{The graph $G'$.}
\label{fig:G'}
\end{figure}
   
If $\widehat{G}$ does not have a trail of length $n+2$, then \poneposs{} winning strategy is to start the game by playing $u \to c$, thus securing a trail of length $n+2$. If \pone{} does so, then \ptwo{} can only get a trail of length at most $n+1$, regardless of whether they begin play on $Q$ or $\widehat{G}$. Thus, \pone{} wins.

If $\widehat{G}$ does have a trail of length $n+2$, then there are several options for \poneposs{} first move.
\begin{enumerate}
\item [(i)] \poneposs{} first move is $x \to y$ for some $xy \in E(\widehat{G})$.\\  If \ptwo{} moves $u \to c$, then since $G$ is cubic and $\widehat{G}$ is defined as above, $\MaxTrail(G'-\{xy,uc\}, y) \leq n+1$, and clearly $\MaxTrail(G'-\{xy,uc\}, c) = n+1$. Thus by Lemma \ref{lem P2 and P1 never meet on a graph}, \ptwo{} wins.

\item [(ii)] \poneposs{} first move is $c \to u$ or $v_{i \pm 1} \to v_{i}$, where $v_i \neq c$.\\
        Since $\deg(u) = 2$ and $\deg(v_i) \leq 2$, \ptwo{} wins by Proposition \ref{lem:p1 loses on degree 1 and 2}. 

\item [(iii)] \poneposs{} first move is $x \to c$ for some $x \in V(G')$.\\ This means the longest trail \pone{} can get is of length $n+2$. \ptwoposs{} winning strategy is to make their first move by playing $u\to w$ (possible no matter which of the three options is \poneposs{} first move), and then continue along a trail of length $n+2$ in $\widehat{G}$. Thus, \ptwo{} secures a trail of length $n+2$ and wins.
\end{enumerate}    

    Thus, \ptwo{} wins {\sc Trail Trap} on $G'$ if and only if $\widehat{G}$ has a trail of length $n+2$, which implies that {\sc Trail Trap \ptwo{}-win} is NP-hard on connected bipartite planar graphs with maximum degree 4.
\end{proof}

\section{Open Problems}
\label{section:conclusion}
There are a number of problems remaining, some of which we mention here. 

\begin{question}
    For which $n\in \mathbb{N}$ is $K_n$ a \pone{}-win graph?
\end{question}
For $n > 9$, determining the winner of {\sc Trail Trap} on $K_n$ requires excessive computation time with our current implementation. Similarly, we are unable to compute the winner on complete bipartite graphs exceeding $16$ vertices (besides those covered by Theorem \ref{prop:complete_bipartite}). However, we pose a conjecture for balanced complete bipartite graphs. Since $K_{3,3}, K_{5,5}$, and $K_{7,7}$ are each \ptwo{}-win, we conjecture the same in general. 

\begin{conjecture}\label{conj:Knn}
   For $n\ge 3$ odd, $K_{n,n}$ is \ptwo{}-win.
\end{conjecture}

In Section 3, we prove that grid graphs $P_m \,\square \, P_n$ are \pone{}-win if $m = 2$ and $n$ is an odd integer at least $5$, but more generally \ptwo{}-win where $m$ and $n$ have the same parity. We conjecture that all larger grid graphs are \ptwo{}-win. 

\begin{conjecture}
    For any $m,n \geq 3$, the grid graph $P_m \mathbin{\square} P_n$ is \ptwo{}-win.
\end{conjecture}

We show above that {\sc Trail Trap \ptwo{}-win} is NP-hard  even on rather restricted classes of graphs. We do not prove NP-completeness because we suspect that the problem is not in NP. Furthermore, we conjecture the following stronger statement is true. 
\begin{conjecture}\label{conj:PSPACE}
    \label{prob:pspace}
    {\sc Trail Trap \ptwo{}-win} is PSPACE-complete.
\end{conjecture}

Since a study of cubic graphs is crucial to our proof that {\sc Trail Trap \ptwo{}-win} is NP-hard, the following is of interest.
\begin{question}\label{ques:hamiltonian}
    Is there a nontrivial restricted class of cubic graphs that are \pone{}-win? For instance, what can be said about {\sc Trail Trap} on Hamiltonian cubic graphs?
\end{question}

Finally, we find that most small graphs are \ptwo{}-win. Indeed, much of our work focuses on finding graph classes where \pone{} wins. As a result, we conjecture the following.
\begin{conjecture}\label{conj:asymptotic}
    The proportion of \pone{}-win graphs on $n$ vertices tends to $0$ as $n\to \infty$. 
\end{conjecture}
The conjecture would be of interest given that we have shown that it is NP-hard to determine whether \ptwo{} wins {\sc Trail Trap} even on restricted families of graphs. 

\section*{Acknowledgments}
Part of this research was conducted at the 2022 Graduate Research Workshop in Combinatorics, which was supported in part by NSF grant 1953985 and a generous award from the Combinatorics Foundation.
M.C.~was supported by the Natural Sciences and Engineering Research Council of Canada (NSERC) Canadian Graduate Scholarship (No. 456422823). 
A.C.~was supported by the Institute for Basic Science (IBS-R029-C1) and by the grant 23-06815M of the Grant Agency of the Czech Republic. 
V.I.C.~was supported by the Slovenian Research and Innovation Agency (ARIS) under the grants Z1-50003, P1-0297, N1-0218, N1-0285, N1-0355, and by the European Union (ERC, KARST, 101071836). 

We thank Eva Arneman, Altea Catanzaro, and Saideh Danison, of the Girls' Angle Math Club, for the early development of {\sc Trail Trap} with mentor R.W. We also thank the reviewers for their helpful comments and suggestions.

\end{document}